\documentclass[]{siamonline181217mod}

\usepackage{amssymb,amsfonts,amsmath}

\usepackage[mathscr]{eucal}

\usepackage{stmaryrd}

\usepackage{amsbsy}

\usepackage{bm}

\usepackage{braket}

\usepackage{multirow}

\usepackage{xfrac}

\usepackage{xparse}

\usepackage{environ}

\usepackage{upquote}

\usepackage{moreverb}

\usepackage{placeins}

\usepackage{enumitem}

\usepackage{tikz}
\usetikzlibrary{calc}
\usepackage{pgfplots}

\usepackage[font=footnotesize,justification=Centering,singlelinecheck=false]{subfig}

\usepackage{cleveref}

\usepackage{setspace}

\definecolor{grey}{RGB}{100,100,100}
\algrenewcommand{\algorithmiccomment}[1]{\hfill\textcolor{grey}{// #1}}

\NewDocumentCommand \qtext {m} {\quad\text{#1}\quad}

\NewDocumentCommand \Real {} {\mathbb{R}}

\NewDocumentCommand \Range { m } { \set{1,2,\dots,#1} }

\NewDocumentCommand \T { O{} m } {\boldsymbol{#1\mathscr{\MakeUppercase{#2}}}}

\NewDocumentCommand \Mx { O{} m } {{\bm{#1\mathbf{\MakeUppercase{#2}}}}} 

\NewDocumentCommand \V { O{} m } {{\bm{#1\mathbf{\MakeLowercase{#2}}}}}

\NewDocumentCommand \X { } {\T{X}}

\NewDocumentCommand \nzX {} {\nnz{\X}}

\NewDocumentCommand \Xk { G{k} } {\Mx{X}_{(#1)}}

\NewDocumentCommand \xe { } {x_{i}}

\NewDocumentCommand \M {} {\T{M}}
\NewDocumentCommand \Mtrue {} {\M_{\rm true}}
\NewDocumentCommand \me { } {m_{i}}

\NewDocumentCommand \FXM {} {F(\X,\M)}

\NewDocumentCommand \fxm {t' G{x} G{m}} {\IfBooleanTF{#1}{\FPD{f}{m}}{f}(#2,#3)}

\NewDocumentCommand \fxme {t'} {\IfBooleanTF{#1}{\fxm'{\xe}{\me}}{\fxm{\xe}{\me}}}

\NewDocumentCommand \we {} {w_i}
\NewDocumentCommand \Y {} {\T{Y}}
\NewDocumentCommand \Ys {} {\T[\tilde]{Y}}

\NewDocumentCommand \Yk { O{k} } {\Mx{Y}_{(#1)}}
\NewDocumentCommand \Yks { O{k} } {\Mx[\tilde]{Y}_{(#1)}} %

\NewDocumentCommand \ye { } {y_{i}}
\NewDocumentCommand \yes { } {\tilde y_i}

\NewDocumentCommand \Ak { G{k} t' t"  } { \Mx{A}_{#1}\IfBooleanTF{#2}{^{\intercal}}{}\IfBooleanTF{#3}{^{\phantom{\intercal}}}{} }
\NewDocumentCommand \EstAk {} {\Mx[\hat]{A}_k}
\NewDocumentCommand \Akj {O{k} G{j}} {\V{a}_{#1}(:,#2)}
\NewDocumentCommand \EstAkj {O{k} G{j}} {\V[\hat]{a}_{#1}(:,\pi(#2))}

\NewDocumentCommand \Akset { } {\set{\Ak }}
\NewDocumentCommand \Bkset { } {\set{\Bk }}
\NewDocumentCommand \Ckset { } {\set{\Ck }}

\NewDocumentCommand \Gk { G{k}  } { \Mx{G}_{#1} }
\NewDocumentCommand \Gks { G{k}  } { \Mx[\tilde]{G}_{#1} } %

\NewDocumentCommand \Gksset { } {\set{\Gks}}

\NewDocumentCommand \gvec {} {\V{g}}
\NewDocumentCommand \gvecs {} {\V[\tilde]{g}}

\NewDocumentCommand \KT { s } {
  \llbracket
  \IfBooleanTF{#1}{\lvec;}{}
  \Ak{1}, \Ak{2}, \dots,  \Ak{d} \rrbracket
}

\NewDocumentCommand \Bk { G{k} s } { \IfBooleanTF{#2}{\Mx[\hat]{B}_{#1}}{\Mx{B}_{#1}} }
\NewDocumentCommand \Ck { G{k} s } { \IfBooleanTF{#2}{\Mx[\hat]{C}_{#1}}{\Mx{C}_{#1}} }

\NewDocumentCommand \Zk { G{k} t' t"} {\Mx{Z}_{#1}\IfBooleanTF{#2}{^{\intercal}}{}%
  \IfBooleanTF{#3}{^{\phantom{\intercal}}}{}}

\NewDocumentCommand \ZkZkt { G{k} } {\Zk{#1}'\Zk{#1}"}

\NewDocumentCommand \SIcnt {} {\tilde{s}_i}
\NewDocumentCommand \NIcnt {} {\tilde{p}_i}
\NewDocumentCommand \ZIcnt {} {\tilde{q}_i}

\NewDocumentCommand \szero {} {s_{\text{\tiny zero}}}
\NewDocumentCommand \sreject {} {s_{\text{\tiny reject}}}

\NewDocumentCommand \pzero {} {p_{\text{\tiny zero}}}
\NewDocumentCommand \FPD { s m m } {
  \IfBooleanTF{#1}
  {\tfrac{\partial #2}{\partial #3}}
  {\frac{\partial #2}{\partial #3}}
}

\DeclareMathOperator{\rank}{rank}
\NewDocumentCommand{\vc}{}{\textsc{vec}}
\NewDocumentCommand{\Exp}{m}{\mathbb{E}[#1]}

\NewDocumentCommand{\nnz}{m}{\text{nnz}(#1)}

\NewDocumentCommand{\LineFor}{m m}{%
  \State\textbf{for} {#1}, \textbf{do} {#2}, \textbf{end}
  }

\NewDocumentCommand \plow {} {\rho_{\rm low}}
\NewDocumentCommand \phigh {} {\rho_{\rm high}}

\newcommand{\tikzmark}[1]{\tikz[overlay,remember picture] \node (#1) {};}

\newcommand*{\AddNote}[5]{%
    \begin{tikzpicture}[overlay, remember picture]
        \draw [decoration={brace,amplitude=0.5em},decorate,grey]
            ($(#3)!(#1.north)!($(#3)-(0,1)$)$) --
            ($(#3)!(#2.south)!($(#3)-(0,1)$)$)
            node [align=center, text width=#4, pos=0.5, anchor=west, right=1em] {//
              #5};
    \end{tikzpicture}
}%

\newcommand{\rundesc}{Each dashed line represents a single run, and the markers signify epochs. The marker is an asterisk  if the true solution was recovered and a dot otherwise. Solid lines represent the median. Dashed black line is the function value estimate for the true solution.}

\newcommand{\lossdesc}{The \emph{same} set of samples is used to estimate the loss across every individual run.}

\NewDocumentCommand{\cpic}{m}{\includegraphics[trim=0 0 15 0, clip, scale=0.4]{fig-chicago-details-#1}}

\NewDocumentCommand{\morecomppic}{m}{%
  \begin{figure}%
    \centering
    \includegraphics[trim=0 15 15 20, scale=0.3]{fig-chicago-semistrat-#1}
    \caption{Component #1 for Chicago crime tensor using semi-stratified sampling.}
    \label{fig:comp-#1}%
  \end{figure}%
}
  
\Crefname{ALC@unique}{Line}{Lines}

\title{Stochastic Gradients for Large-Scale Tensor Decomposition%
  \thanks{This work was supported by the Laboratory Directed Research
    and Development Program at Sandia National Laboratories and
    the U.S. Department of Energy, Office of Science,
    Office of Advanced Scientific Computing Research (ASCR)
    Applied Mathematics Program.
    Sandia National Laboratories is a multimission laboratory
    managed and operated by National Technology and Engineering
    Solutions of Sandia, LLC., a wholly owned subsidiary of Honeywell
    International, Inc., for the U.S. Department of Energy's National
    Nuclear Security Administration under contract DE-NA-0003525.
    This paper describes objective technical results and analysis. Any
    subjective views or opinions that might be expressed in the paper
    do not necessarily represent the views of the U.S. Department of
    Energy or the United States Government.
    Work by D.H. was also supported in part by NSF Grant IIS 1838179.}}

\author{%
  Tamara~G.~Kolda%
  \thanks{Sandia National Laboratories, Livermore, CA (\email{tgkolda@sandia.gov})}
  \and David Hong%
  \thanks{University of Michigan, Ann Arbor, MI (\email{dahong@umich.edu})}
 }

\headers{Stochastic Gradients for Tensor Decomposition}
{Tamara G. Kolda and David Hong}

\begin{document}

\maketitle

\begin{abstract}
  Tensor decomposition is a well-known tool for multiway data analysis.
  This work proposes using stochastic gradients for efficient generalized canonical polyadic (GCP) tensor decomposition of large-scale tensors.
  GCP tensor decomposition is a recently proposed version of tensor decomposition that allows for a variety of loss functions such as Bernoulli loss for binary data or Huber loss for robust estimation.
  The stochastic gradient is formed from randomly sampled elements of the tensor and is efficient because it can be computed using the sparse matricized-tensor-times-Khatri-Rao product (MTTKRP) tensor kernel.
  For dense tensors, we simply use uniform sampling.
  For sparse tensors, we propose two types of stratified sampling that give precedence to sampling nonzeros.
  Numerical results demonstrate the advantages of the proposed approach and its scalability to large-scale problems.
\end{abstract}


\section{Introduction}
\label{sec:introduction}

Tensor decomposition is the higher-order analogue of matrix
decomposition and is becoming  an everyday tool for data analysis.
It can be used for unsupervised learning, dimension reduction, tensor completion, feature extraction in supervised machine learning, data visualization, and more.
For a given $d$-way data tensor $\X$ of size $n_1 \times n_2 \times \cdots \times n_d$,
the goal is to find a \emph{low-rank} approximation $\M$, i.e.,
\begin{equation}\label{eq:Mdef}
  \X \approx \M \qtext{where}
  \M = \sum_{j=1}^r \Akj[1] \circ \Akj[2] \circ \cdots \circ \Akj[d].
\end{equation}
The low-rank structure of $\M$ reveals patterns within the data, as defined by the \emph{factor matrices}.
The $k$th factor matrix of size $n_k \times r$ is denoted $\Ak$.
The $j$th \emph{factor} (column) in mode $k$ is denoted $\Akj$.
Each \emph{component} of $\M$ is a $d$-way outer product (denoted by $\circ$) of $d$ factors, forming a \emph{rank-one} tensor.
We say $\rank(\M) \leq r$
because $\M$ can be written as the sum of $r$ rank-one tensors.
We define $n \equiv (\prod_{k=1}^d n_k)^{1/d}$ and  assume $r \ll n^d$.
The storage for $\X$ is $n^d$ if it is dense and $O(\nnz{\X})$ if it is sparse.
The storage for $\M$ is $O(r  \sum_{k=1}^d n_k)$, which is usually much
smaller than that for $\X$.

The standard canonical polyadic or CANDECOMP/PARAFAC (CP) tensor decomposition seeks the best low-rank approximation with respect to sum of squared errors \cite{CaCh70,Ha70}.
The generalized CP (GCP) tensor decomposition is a novel approach that allows for an arbitrary elementwise loss function that is summed across all tensor entries \cite{HoKoDu20}; i.e., the user provides a scalar loss function $f: \Real\times\Real \rightarrow \Real$ resulting in the optimization problem
\begin{equation}
  \label{eq:gcp}
  \text{minimize } \FXM \equiv \sum_{i_1=1}^{n_1} \cdots \sum _{i_d=1}^{n_d} \fxm{x_{i_1\dots i_d}}{m_{i_1\dots i_d}} \qtext{subject to} \rank(\M) \leq r.
\end{equation}
GCP is useful in situations where non-standard choices of the scalar loss function $f$ may be appropriate.
For instance, if the entries of $\X$ are binary values in $\set{0,1}$, we may use the Bernoulli loss, i.e.,
$\fxm = \log(m+1)-x\log(m)$
where $m$ represents the \emph{odds} of observing a one.
For positive data that has a Gamma distribution, we may
use $\fxm=x/m + \log(m)$.
Standard CP tensor decomposition corresponds to $\fxm = (x {-} m)^2$, and
Poisson tensor decomposition \cite{WeWe01,ChKo12} to $\fxm = m - x \log m$.
See \cite{HoKoDu20} for full details.

In this paper, we consider the problem of fitting GCP for large-scale tensors.
The GCP gradient involves a sequence of $d$
matricized-tensor times Khatri-Rao products ({MTTKRPs})
with a \emph{dense} tensor of size $n^d$ that costs $O(rn^d)$ operations, even when $\X$ is sparse.%
\footnote{Both standard and Poisson CP have special structure that allows the gradient to be formed implicitly, as described in \cref{sec:special-cases} so that sparsity of $\X$ can be exploited.}
Thus, the computational and
storage costs of computing the gradient may be prohibitive.
Our solution is a flexible framework for \emph{stochastic} gradient computation for GCP,
which can be used with stochastic gradient descent (SGD) or popular variants such as  Adam \cite{KiBa15}.
Our stochastic framework replaces the dense tensor with a (stochastic) sparse tensor that equals the dense tensor in expectation.
If the stochastic gradient samples $s \ll n^d$ tensor entries,
the resulting \emph{sparse} MTTKRP
reduces the gradient cost to $O(sr)$ flops and $O(s)$ intermediate storage%
, albeit with a potentially significant
sacrifice in accuracy due to stochasticity.

A distinguishing feature of our framework is stratified sampling.
If $\X$ is sparse, uniform sampling of the indices rarely samples nonzeros even though they are extremely important to
the minimization.
Stratified sampling randomly selects nonzeros disproportionately,
which can reduce the variance of the stochastic gradient, accelerating convergence.
Since stratified sampling can itself be expensive, we also introduce a semi-stratified sampling approach to
further improve efficiency.
Note that simply ignoring the zero entries, as is done in recommender systems where zeros indicate missing data, yields a fundamentally different problem.
\Cref{sec:related-work} surveys related work,
\cref{sec:notation-background} provides notation and background,
and \cref{sec:stoch-gcp-grad} discusses the framework.

We present computational experiments on both artificial and real-world problems
in \cref{sec:experimental-results}, demonstrating the reliability and efficiency
of using stochastic gradients.
Experimentally, we can set the number of gradient samples to be the sum of the tensor dimensions, i.e., $\sum_{k=1}^d n_k$.
For dense problems, stochastic optimization can be an order of magnitude faster than non-stochastic optimization.
For sparse problems, stochastic optimization enables computing GCP for sparse tensors that would otherwise be intractable.
Additionally, stratified and semi-stratified sampling of zeros and nonzeros have a clear advantage over uniform sampling in this setting.
We also compare to CP-APR \cite{ChKo12,HaPlKo15} on real-world Chicago crime data. In our experiments, the stochastic methods are
faster and find equivalent solutions.

\section{Related Work}
\label{sec:related-work}

This work is a follow-on to the introduction of GCP by Hong, Kolda, and Duersch \cite{HoKoDu20}.
They only considered small-scale problems that could be treated as dense.
The goal of this work is to tackle larger scale problems via stochastic gradient methods.

A major motivation for the present work was that of Acar et al.~\cite{AcDuKoMo11} which focused on standard CP factorization
for tensors with irretrievably missing data.
They considered cases where the vast majority of the data is missing,
i.e., \emph{scarce} tensors in the terminology of \cite{HoKoDu20},
and observed that a scarce data tensor leads to a sparse MTTKRP in computing the gradient.
That observation proves to be important here because our sampled data tensors are scarce.

In the specific case of third-order ($d=3$) and one sample per iteration ($s=1$),
Beutel et al.~\cite{BeTaKuFa14} have proposed SGD for standard CP tensor decomposition.
Since a single tensor element only updates one row in each of the three factor matrices, they are able to parallelize updates
that do not operate on the same rows.
They do not consider sampling strategies or the MTTKRP structure in the gradient.

Vervliet and De Lathauwer \cite{VeLa16} pursue an alternative form of sampling for standard CP:
rather than sampling \emph{elements}, they sample \emph{indices from each mode} of the tensor.
The resulting random subtensor is used in either one outer iteration of alternating least squares (ALS) or
one iteration of Gauss-Newton, with careful attention to step sizes as
is needed in stochastic optimization.
The updates only modify the rows of the factor matrices corresponding to the subsampled indices.
This block sampling strategy is not readily amenable to sparse data tensors,
which are a focus of our work.

For alternating least squares (CP-ALS), randomized matrix sketching %
can be used to solve the linear least squares problems in the dense \cite{BaBaKo18} and sparse \cite{ChPePeLi16} cases; this approach does not directly apply to GCP because its subproblems are not least squares problems.
Also in the vein of sketching, random projections can be used to build a sketch (or multiple sketches) of the complete tensor \cite{PaFaSi12,SiPaFa14,ZhCiXi14},
including the special case of symmetric factorizations \cite{WaTuSmAn15,SoWoZh16}.

In the streaming context, one dimension is typically treated as time with a new slice arriving at each time step
and can be updated via SGD; see, e.g., \cite{MaMaGi15}.
Maehara, Hayashi, and Kawarabayashi~\cite{MaHaKa16}
factorize tensors that arise as sums of streaming tensor samples.
Similar to our approach,
they apply SGD on the samples to update the factorization of the sum.
In contrast to our focus on effective sampling, they primarily consider applications where a stream of samples
is already given.
There are also other works on updates for streaming tensors that are not based on SGD \cite{NiSi09, VaVeLa17, MaYaWa18, GuPaPa18}.

In recommender systems and tensor completion,
it is typically assumed that most data is missing,
resulting in scarce tensors (sometimes incorrectly conflated with sparse tensors).
Smith, Park, and Karypis~\cite{SmPaKa16} consider SGD for scarce tensors
and only sample nonzeros as they correspond to observed entries.
For symmetric tensors with missing data,
Ge et al.~\cite{GeHuJiYu15} prove that SGD
finds a local minimizer rather than a
saddle point.

We close by highlighting that SGD has already had broad success
in matrix factorization, primarily for
recommender systems where zeros are treated as missing data (scarce matrices); see, e.g., \cite{KoBeVo09} and references therein.
However, we omit a full discussion of stochastic gradients in the matrix setting because the field is vast.
To the best of our knowledge, nothing like the stratified and semi-stratified sampling proposed in this work has
been proposed in the matrix regime.
Gemulla et al.~\cite{GeNiHaSi11} propose a stratified version of SGD
as a technique to partition the data into independent segments that can be processed in parallel,
but their use of stratification is different than our proposal.

\section{Notation and background}
\label{sec:notation-background}

We provide some context
for the remainder of the paper.

\subsection{Probability notation and sampling background}
\label{sec:prob-stat-notation}

In this paper, we use a tilde to indicate random variables and instances thereof, e.g., $\tilde x$.
The \emph{expected value} is denoted $\Exp{\tilde x}$.
If $\tilde x \in \set{v_1, v_2, \dots v_m}$ is a discrete random variable that samples $v_i$ with probability $p_i$, then
\begin{displaymath}
  \Exp{\tilde x} \equiv \sum_{i=1}^m p_i \, v_i.
\end{displaymath}
The random sample is
\emph{uniform} if  each element has equal probability;
i.e., $p_i = 1/m$ for $i=1,2,\dots,m$.
Sampling \emph{with replacement} means that the same $v_i$ can be sampled more than once,
i.e., that subsequent samples are independent identically distributed (i.i.d.) draws.
If we let $\SIcnt$ denote the number of times that $v_i$ is sampled from a uniform distribution over $s$ draws with replacement, then
$\Exp{\SIcnt} = s/m$
since there are $s$ independent draws and each draw selects sample $v_i$ with probability $1/m$.

\subsection{Tensor notation}
\label{sec:tensors-notation}

For a $d$-way data tensor $\X$ of size $n_1 \times n_2 \times \cdots \times n_d$,
we define
\begin{displaymath}
  n = \sqrt[d]{\prod_{k=1}^d n_k}
  \qtext{and}
  \bar n = \frac{1}{d} \sum_{k=1}^d n_k.
\end{displaymath}
These quantities provide convenient measures of how large the tensor is. We treat $d$ as constant for big-O notation.
Throughout, we let $i = (i_1,i_2,\dots,i_d)$ where $i_k \in \Range{n_k}$ and refer to this as a \emph{multi-index} or simply an \emph{index}.
Every multi-index has a corresponding  \emph{linear index} between 1 and $n^d$ \cite{KoBa09}.
We let $\Omega$ denote the set of all tensor indices, so  necessarily $|\Omega| = n^d$.
Excepting the discussion of missing data in \cref{sec:weight-formulations}, we assume all tensor entries are known.
We let $\nzX$ denote the number of nonzeros in $\X$ and say that $\X$ is sparse if $\nzX \ll  n^d$.

If $\X$ is dense, then the storage is $n^d$.
If $\X$ is sparse, then only the nonzeros and their indices are stored, so the storage is $(d+1)\,\nzX$ using a coordinate format \cite{BaKo07}.
The total storage for a rank-$r$ approximation $\M$ as in \cref{eq:Mdef} is the storage for the factor matrices, i.e., $\bar n d r$, and is usually significantly less than the storage for $\X$

\subsection{MTTKRP background}
\label{sec:background-mttkrp}

Given a tensor $\Y$ of size $n_1 \times n_2 \times \cdots \times n_d$ and factor matrices $\Ak$ (from a low-rank $\M$) of size $n_k \times r$ for $k=1,\dots,d$,
the MTTKRP in mode $k$ is
\begin{equation}\label{eq:mttkrp}
  \Yk \underbrace{(\Ak{d} \odot \cdots \odot \Ak{k+1} \odot \Ak{k-1} \odot \cdots \odot \Ak{1})}_{\Zk}.
\end{equation}
The matrix $\Zk$ is of size $n^d/n_k \times r$ and is the Khatri-Rao product (denoted by $\odot$)
of all the factor matrices except $\Ak$, and
the matrix $\Yk$ is of size $n_k \times n^d/n_k$ and denotes the mode-$k$ unfolding of $\Y$.

Much work has gone into efficient computation of MTTKRP.
Bader and Kolda \cite{BaKo07} consider both dense and sparse $\Y$ tensors, showing
that the cost is $O(rn^d)$ for dense $\Y$ and $O(r\,\nnz{\Y})$ for sparse $\Y$.
Phan, Tichavsky, and Cichocki \cite{PhTiCi13} propose methods to reuse partial computations
when computing the MTTKRP for all $d$ modes in sequence.
Much recent work has focused on more efficient representations of
sparse tensors and parallel MTTKRP computations
\cite{SmRaSiKa15,KaUc15,LiChPeSu17,PhKo19}.
There is also continued work on improving the efficiency of dense MTTKRP calculations \cite{HaBaJiTo18,BaKnRo18}.

\section{Stochastic GCP Gradient}
\label{sec:stoch-gcp-grad}

GCP tensor decomposition generalizes CP tensor decomposition, minimizing the GCP loss function in \cref{eq:gcp}.
From \cite[Theorem~3]{HoKoDu20}, the GCP gradient is calculated via a sequence of MTTKRP computations.
Key to this calculation is the \emph{elementwise partial gradient tensor}
$\Y$ that is the same size as $\X$ and is defined as
\begin{equation}
  \label{eq:Y}
  \ye =  \fxme' .
\end{equation}
That is, $\Y$ is the tensor of first derivatives of $f$ evaluated elementwise w.r.t.\@ $\X$ and $\M$.
Whether the data tensor $\X$ is dense or sparse, the elementwise partial gradient tensor $\Y$ is dense almost everywhere.
The partial derivative, denoted $\Gk$, of the objective $F$ in \cref{eq:gcp}  w.r.t.\@ $\Ak$  is
\begin{equation}
  \label{eq:Gk}
  \Gk  = \Yk \Zk,
\end{equation}
where $\Zk$ is the Khatri-Rao product defined in \cref{eq:mttkrp}.
Because $\Y$ is dense (even when $\X$ is sparse),
$\Y$ requires $n^d$ storage and
the cost of computing MTTKRP for gradients $\Gk$ is $O(rn^d)$.
Such costs may be prohibitive.
For instance, if $n=1000$ and $d=4$, then $\Y$ would require 8~TB of storage.
There are a couple of special cases (see \cref{sec:special-cases})
where dense computation can be avoided because $\Y$ is formed \emph{implicitly},
notably for standard CP \cite{BaKo07},
but this is not the case for general loss functions and so motivates stochastic approaches.

The stochastic approach is based on the following observation:
$\Y$ becomes \emph{sparse} when the data tensor $\X$ is \emph{scarce}
(i.e., most of the elements are missing).
This is because,
by \cite[Theorem~3]{HoKoDu20},
the gradient in this case is the same as in \cref{eq:Gk} except that $\Y$ is defined as
\begin{displaymath}
  \ye =
  \begin{cases}
    \fxme' & \text{if $\xe$ known}, \\
    0 & \text{if $\xe$ missing}.
  \end{cases}
\end{displaymath}
If $\Y$ is sparse, then the MTTKRP in \cref{eq:Gk} can be computed in time $O(r \, \nnz{\Y})$.

We develop a stochastic gradient whose general form is
\begin{equation}
  \label{eq:Gks}
  \Gks = \Yks \Zk %
  \qtext{where} \Exp{\Ys} = \Y
  \qtext{ and }  \nnz{\Ys} \leq s \ll n^d.
\end{equation}
By linearity of expectation, $\Exp{\Ys} = \Y$ implies $\Exp{\Gks} = \Gk$.
Making $\Ys$ sparse unlocks efficient sparse computation of the MTTKRP.
Consequently, storage for $\Ys$ is $O(s)$ and the cost to compute the gradients is
$O(rs)$, a reduction of  roughly $n^d / s$ compared to computing the full gradient, albeit
at the cost of lower accuracy.

In the remainder of this section, we discuss the pros and cons of several different choices for $\Ys$ based on
uniform sampling, stratified sampling, and semi-stratified sampling.

\subsection{Uniform Sampling}
\label{sec:uniform-sampling}
To create a random instance of $\Ys$,
we sample $s$ indices uniformly with replacement.
The number of times that $i$ is sampled is denoted as $\SIcnt$, so $\sum_{i\in\Omega} \SIcnt = s$.
Using this sampling strategy, in expectation we have
\begin{displaymath}
  \Exp{\SIcnt} = \frac{s}{n^d}
  \qtext{for all} i \in \Omega.
\end{displaymath}
Note that the $\SIcnt$ values are not stored as a dense object but rather as just a list of the samples or in some other sparse data structure.
We then define the stochastic tensor $\Ys$ to be
\begin{displaymath}
  \yes = \SIcnt \, \frac{n^d}{s}  \, \ye
  \qtext{for all} i \in \Omega.
\end{displaymath}
The stochastic $\Ys$ is sparse because at most $s$ entries are nonzero (since
at most $s$ values of $\SIcnt$ are nonzero).
Clearly $\Exp{\Ys} = \Y$ since
$\Exp{\yes} = \Exp{\SIcnt} \, \frac{n^d}{s} \, \ye = \ye$
for all $i \in \Omega$.
A similar argument applies for sampling \emph{without} replacement,
but there is no practical difference when $s \ll n^d$.

The simplest way to generate a random index $i$ is to generate $d$ random mode indices $i_k$:
\begin{displaymath}
  i_k = \texttt{randi}(n_k)
  \qtext{for} k \in \Range{d}.
\end{displaymath}
Here, $\texttt{randi}(m)$ indicates selecting a random integer between 1 and $m$.
This requires $d$ random numbers per index $i$.
Alternatively, we can generate a random \emph{linear} index via $\texttt{randi}(n^d)$ and then
convert it to a tensor multi-index.%
\footnote{Generating $d$ separate entries helps prevent overflow if $n^d > 2^{64}$,
  i.e., the total number of entries in the tensor is more than the size of the largest unsigned $64$-bit integer (\texttt{uint64} in MATLAB).
  Larger integers can be generated by instead using extended-precision arithmetic
  (e.g., \texttt{int} in Python).}
Either way is $O(1)$ work (we treat $d$ as a constant).

The procedure is presented in \cref{alg:stoc-grad-uniform}.
The function \texttt{MTTKRP} corresponds to
\cref{eq:mttkrp}, and this can be computed efficiently because
specialized sparse implementations exist as discussed in
\cref{sec:background-mttkrp}.
The model entries $\me$ and elementwise partial gradient tensor entries $\ye$ are only computed for the $s$ randomly selected indices.
The most expensive operations are computing the model entries $\me$ and the corresponding MTTKRP calculations.
In the implementation, we are able to share some intermediate computations across these two steps.

\begin{algorithm}
  \caption{Stochastic Gradient with Uniform Sampling}
  \label{alg:stoc-grad-uniform}
  \begin{algorithmic}[1]\footnotesize
    \Function{StocGrad}{$\X, \Akset, s$}
    \State{\label{line:sgu-omega}$\omega \gets \prod_k n_k$} \Comment{$\omega = $ \# entries in $\X$}
    \State $\Ys \gets 0$ \Comment{initialize $\Ys$ to all-zero \emph{sparse} tensor}
    \For{$c = 1,2,\dots,s$} \Comment{loop to sample $s \ll \omega$  indices for $\Ys$}
    \LineFor{\label{line:sgu-index}$k = 1,2,\dots,d$}{$i_k \gets \texttt{randi}(n_k)$} \Comment{sample random index $i \equiv (i_1,i_2,\dots,i_d)$}
    \State $\me \gets \sum_{j=1}^r \prod_{k=1}^d a_k(i_k,j)$ \Comment{compute $\me$ at sampled index}
    \State\label{line:sgu-yes}$\yes \gets \yes +  (\omega/s) \,  g(\xe,\me)$ \Comment{compute $\yes$ at sampled index, $g \equiv \partial f / \partial m$}
    \EndFor{\label{line:sgu-endfor}}
    \LineFor{$k=1,2,\dots,d$}{$\Gks \gets \mathtt{MTTKRP}(\Ys, \Ak, k)$} \Comment{use stochastic sparse $\Ys$ to compute $\Gks$}
    \State \Return $\Gksset$
    \EndFunction
  \end{algorithmic}
\end{algorithm}

\subsection{Stratified Sampling}
\label{sec:stratified-sampling-1}

Uniform sampling may not be appropriate for sparse tensors since nonzeros will rarely be sampled
(each draw gets a nonzero with probability $\nzX / n^d$).
Intuitively, however, we expect nonzeros in a sparse tensor to be important to the factorization.

Non-uniform sampling has been used in other stochastic gradient contexts.
Needell et al.~\cite{NeSrWa15} analyze SGD and argue that
biased sampling toward functionals with larger Lipschitz constants can improve performance.
Gopal~\cite{Go16} similarly biases sampling toward functionals with larger gradients
to reduce the variance of the stochastic gradients,
and similar ideas appear in \cite{ZhZh14a,ZhZh14b}.
This idea has many potential applications in machine learning,
and roughly translates in our case
to up-sampling a tensor entry if the corresponding elementwise loss function $f_i \equiv \fxme$
has a higher Lipschitz constant for fixed $\xe$.
Consider GCP with the Bernoulli odds loss \cite{HoKoDu20}:
$\fxm = \log(m+1) - x \log m$ where $x \in \set{0,1}$ and $m > 0$ corresponds to the odds of $x=1$.
For $x=0$, $|\fxm'{0}{m}| = 1/(m+1) \leq 1$;
for $x=1$, $|\fxm'{1}{m}| = 1/(m^2 + m)$ is unbounded as $m \downarrow 0$.
Thus, the elementwise loss functions for nonzeros are not Lipschitz,
their gradients can be very large,
and sampling them more often by sampling zeros and nonzeros separately
can reduce the variance of the stochastic gradients.

Consider a generic partition of $\Omega$ into $p$ partitions $\Omega_1, \Omega_2, \dots, \Omega_p$.
If we partition into zeros and nonzeros, then $p=2$.
Let $s_\ell > 0$ be the number of samples from partition $\Omega_{\ell}$
so that the total number of samples is $s = \sum_{\ell} s_{\ell}$.
Within each partition, we sample uniformly with replacement.
Thus, the expected number of times index $i$ is sampled depends on its partition:
\begin{displaymath}
  \Exp{\SIcnt} = \frac{s_{\ell}} {|\Omega_{\ell}|}
  \qtext{where} i \in \Omega_{\ell}
    \qtext{for all} i \in \Omega.
\end{displaymath}
We then define the stochastic tensor $\Ys$ to be
\begin{displaymath}
  \yes = \SIcnt \, \frac{|\Omega_{\ell}|}{s_{\ell}}  \, \ye
  \qtext{where} i \in \Omega_{\ell}
  \qtext{for all} i \in \Omega.
\end{displaymath}
Once again,  at most $s$ entries of $\Ys$ are nonzero, and $\Exp{\Ys} = \Y$ since
$\Exp{\yes} = \Exp{\SIcnt} \, \frac{|\Omega_{\ell}|}{s_{\ell}}  \, \ye = \ye$
{for all} $i \in \Omega$.
Interestingly, the weight (${|\Omega_{\ell}|}/{s_{\ell}}$) for index $i$ is completely
independent of the overall sample or tensor size --- it depends only
on the partition size and the number of samples for that partition.

\Cref{alg:stoc-grad-stratified} presents the stratified sampling procedure for sparse tensors,
where we sample $p$ nonzeros and $q$ zeros.
Sampling nonzeros is straightforward because they are stored as a list in coordinate format.
However, sampling zeros requires sampling a random index as was done for uniform sampling
and then \emph{rejecting} the sample if it is actually a nonzero.

\begin{algorithm}
  \caption{Stochastic Gradient with Stratified Sampling for Sparse Tensor}
  \label{alg:stoc-grad-stratified}
  \begin{algorithmic}[1]\footnotesize
    \Function{StocGrad}{$\X, \Akset, p, q$}
    \State $\eta \gets \nzX$ \Comment{$\eta = $ \# of nonzeros in $\X$}
    \State $\zeta \gets \prod_k n_k - \nzX$ \Comment{$\zeta = $ \# of zeros in $\X$}
    \State $\Ys \gets 0$
    \For{$c = 1,2,\dots,p$} \Comment{loop to sample $p \ll (\eta+\zeta)$  nonzero indices for $\Ys$}
    \State $\xi \gets \texttt{randi}(\eta)$ \Comment{sample random nonzero index in $\set{1,\dots,\eta}$}
    \State $i \gets $ index of $\xi$th nonzero \Comment{extract corresponding tensor index}
    \State $\me \gets \sum_{j=1}^r \prod_{k=1}^d a_k(i_k,j)$ \Comment{compute $\me$ at sampled index}
    \State\label{line:sgs-yes-a}$\yes \gets \yes +  (\eta/p) \,  g(\xe,\me)$ \Comment{compute $\yes$ at sampled index, $g \equiv \partial f / \partial m$}
    \EndFor
    \State $c \gets 0$
    \While{$c < q$} \Comment{loop to sample $q \ll (\eta + \zeta)$  zero indices for $\Ys$}
    \LineFor{$k = 1,2,\dots,d$}{$i_k \gets \texttt{randi}(n_k)$} \Comment{sample random index $i \equiv (i_1,i_2,\dots,i_d)$}
    \If{$x_i \neq 0$} \Comment{check against list of nonzeros}
    \State reject sample
    \Else
    \State $c \gets c+1$ \Comment{increment count of accepted zero samples}
    \State $\me \gets \sum_{j=1}^r \prod_{k=1}^d a_k(i_k,j)$ \Comment{compute $\me$ at sampled index}
    \State\label{line:sgs-yes-b}$\yes \gets \yes +  (\zeta/q) \,  g(\xe,\me)$ \Comment{compute $\yes$ at sampled index, $g \equiv \partial f / \partial m$}
    \EndIf
    \EndWhile
    \LineFor{$k=1,2,\dots,d$}{$\Gks \gets \mathtt{MTTKRP}(\Ys, \Ak, k)$} \Comment{use stochastic sparse $\Ys$ to compute $\Gks$}
    \State \Return $\Gksset$
    \EndFunction
  \end{algorithmic}
\end{algorithm}

We can estimate how many random multi-indices are needed to obtain $q$ valid zero samples.
Let
$\eta = \nzX$ and $\zeta = n^d - \eta$.
We expect a very small proportion ($\eta / n^d$) of indices to be rejected.
The number of samples needed to produce an average of $q$ zero
indices is
\begin{equation}\label{eq:oversample-rate}
  \frac{q}{1 - \eta/n^d} = \frac{n^d}{\zeta}\, q\approx q.
\end{equation}
This is only on average, so we oversample to get sufficiently many zeros
with high probability.
Namely, we sample $\rho \, (n^d/\zeta) \, q$ indices where $\rho > 1$ is the oversample rate.
A default of $\rho = 1.1$ is justified in \cref{sec:determ-overs-rate}.

Since we have to check against the entire list of nonzeros for each sample,
the cost of the rejection sampling for zeros can be significant.
To achieve a speed
that does not dominate the other costs,
our MATLAB implementation uses various efficiencies
such as conversion from multi-indices to linear indices,
presorting the list of nonzero linear indices, and using the hidden builtin function \texttt{\_ismemberhelper}.
If $n^d \geq 2^{64}$, conversion to linear indices so that we can use these efficiencies
is not possible.
Beyond MATLAB, efficiency can be achieved using extended precision,  
hash tables, etc.
As an alternative to \emph{specialized} implementations for efficient rejection sampling (that are potentially language and architecture dependent),
we propose a \emph{specialized} algorithm that entirely avoids rejection sampling in the next subsection.

\subsection{Semi-Stratified Sampling}
\label{sec:semi-strat-sampl}

To avoid rejection sampling entirely, we propose a variant that we call \emph{semi-stratified} sampling.
We sample zeros incorrectly but then correct for it when sampling nonzeros.
Specifically, we sample without rejection to obtain ``zeros'', knowing
that a small proportion may actually be nonzeros.
We still sample nonzeros explicitly but now add a correction to account for the possibility that they were also wrongly sampled as ``zeros.''

Namely, we sample $p$ nonzeros and $q$ ``zeros''  from $\Omega$ (the entire set of indices).
Let $\NIcnt$ be the number of times that index $i$ was sampled as a nonzero where $\sum_i \NIcnt = p$. Clearly, $\NIcnt = 0$ if $\xe = 0$.
Let $\ZIcnt$ be the number of times that index $i$ is sampled as a ``zero'' from the full set of possible indices where $\sum_i \ZIcnt = q$.
It is possible that some nonzeros are sampled, i.e., we can have $\ZIcnt > 0$ when $\xe \neq 0$.
Using these counts and recalling $\eta = \nzX$, we define $\Ys$ as
\begin{displaymath}
  \yes = \NIcnt \, \frac{\eta}{p} \, ( \ye - c_i) + \ZIcnt \, \frac{n^d}{q} c_i
  \qtext{where} c_i \equiv \fxm'{0}{\me}.
\end{displaymath}
This still satisfies $\Exp{\yes} = \ye$. For  $i$ such that $\xe = 0$, we have $\NIcnt = 0$ and $c_i = \ye$, so
\begin{displaymath}
  \Exp{\yes} = \Exp{\ZIcnt} \, \frac{n^d}{q} c_i = \frac{q}{n^d}  \, \frac{n^d}{q} \ye = \ye.
\end{displaymath}
For $i$ such that $\xe \neq 0$, we have
\begin{displaymath}
  \Exp{\yes} = \Exp{\NIcnt} \, \frac{\eta}{p} \, ( \ye - c_i) + \Exp{\ZIcnt} \, \frac{n^d}{q} c_i
  = \frac{p}{\eta} \, \frac{\eta}{p} \, ( \ye - c_i) + \frac{q}{n^d} \, \frac{n^d}{q} c_i
  = (\ye - c_i) + c_i = \ye.
\end{displaymath}

The procedure is implemented in \cref{alg:stoc-grad-semistrat}.
The procedure for the nonzero samples is identical to that in \cref{alg:stoc-grad-stratified} except for the adjustment term $-g(0,\me)$ to ensure that the
expectations are correct.
The procedure for the ``zeros'' differs because it samples over the entire index space and does not reject nonzeros; it assumes that $\xe=0$ in computing $\ye$.

\begin{algorithm}
  \caption{Stochastic Gradient with Semi-Stratified Sampling for Sparse Tensor}
  \label{alg:stoc-grad-semistrat}
  \begin{algorithmic}[1]\footnotesize
    \Function{StocGrad}{$\X, \Akset, p, q$}
    \State $\eta \gets \nzX$ \Comment{\# of nonzeros in $\X$}
    \State $\omega \gets \prod_k n_k$ \Comment{\# entries in $\X$}
    \State $\Ys \gets 0$ %
    \For{$c = 1,2,\dots,p$} \Comment{loop to sample $p \ll \omega$  indices for $\Ys$}
    \State $\xi \gets \texttt{randi}(\eta)$ \Comment{sample random nonzero index in $\set{1,\dots,\eta}$}
    \State $i \gets $ index of $\xi$th nonzero \Comment{extract corresponding tensor index}
    \State $\me \gets \sum_{j=1}^r \prod_{k=1}^d a_k(i_k,j)$ \Comment{compute $\me$ at sampled index}
    \State $\yes \gets \yes +  (\eta/p) \,  [g(\xe,\me) - g(0,\me)]$ \Comment{compute $\yes$ at sampled index with \emph{correction} term}
    \EndFor
    \For{$c = 1,2,\dots,q$} \Comment{loop to sample $q \ll \omega$  ``zero'' indices for $\Ys$}
    \LineFor{$k = 1,2,\dots,d$}{$i_k \gets \texttt{randi}(n_k)$} \Comment{sample index $i \equiv (i_1,i_2,\dots,i_d)$}
    \State $\me \gets \sum_{j=1}^r \prod_{k=1}^d a_k(i_k,j)$ \Comment{compute $\me$ at sampled index}
    \State $\yes \gets \yes +  (\omega/q) \,  g(0,\me)$ \Comment{compute $\yes$ at sampled index, \emph{assuming} $\xe=0$}
    \EndFor
    \LineFor{$k=1,2,\dots,d$}{$\Gks \gets \mathtt{MTTKRP}(\Ys, \Ak, k)$} \Comment{use sparse $\Ys$ to compute $\Gks$}
    \State \Return $\Gksset$
    \EndFunction
  \end{algorithmic}
\end{algorithm}

\subsection{Adapting to Weighted Formulations}
\label{sec:weight-formulations}
Any of the above sampling methods can be easily adapted to the weighted version of GCP
with the \emph{weighted} loss function
\begin{equation}
  \label{eq:wgcp}
  \text{minimize } \FXM \equiv \sum_{i} \we \, \fxme \qtext{subject to} \rank(\M) \leq r.
\end{equation}
In estimating the gradient, the only change in the methods is to replace $\ye$ in \cref{eq:Y} with
\begin{displaymath}
  \ye = \we \, \fxme'.
\end{displaymath}
We do not explicitly study the weighted formulation in this work, but we briefly mention two scenarios where such methods may be useful.

\paragraph{Adapting to Weighted Formulations}

In the context of recommender systems, zeros are generally treated as missing data.
The data is assumed to be missing at random (MAR), meaning that the probability of a data item being present does not depend on its value.
However, it has been argued in the matrix case that this may be a flawed assumption~\cite{MaZeRoSl12}.
Instead, we may consider including the zero terms but \emph{down-weighting} them by using a weighted scheme, e.g., $\we=1$ for $\xe \neq 0$ and $\we=0.1$ for $\xe=0$. Many large-scale tensor applications have a recommender system flavor, where down-weighting of the zero entries may be appropriate.
More generally, down-weighting can be appropriate for entries that are less reliable or noisier, as has been done in various matrix factorization applications \cite{CoHo77,JaHoBoGrSm04,TaMaZu05,YuTo04}.
See also \cite{SrJa03,UdHoZaBo16} for work on computing weighted matrix factorizations and \cite{HoFeBa18} for a recent analysis of matrix factorization weights when columns have heterogeneous noise levels.

\paragraph{Weighted Formulations for Missing Data}

Conversely, if some portion of data is missing, there are various strategies that can be used to avoid sampling missing elements.
However, this can also be handled easily by setting the weights of missing entries to be zero.
Ideally, elements with a weight of zero should be avoided during sampling.

\FloatBarrier

\section{Experimental Results}
\label{sec:experimental-results}

All experiments were run using MATLAB (Version 2018a) on a Dual Socket Intel E5-2683v3 2.00GHz CPU with 256~GB memory. %
The methods are implemented in \texttt{gcp\_opt} in the Tensor Toolbox for MATLAB \cite{TensorToolbox}.

\subsection{Stochastic optimization algorithm}
\label{sec:stoch-optim-algor}

The stochastic gradients can be used with any number of stochastic optimization methods.
We use the popular Adam \cite{KiBa15} method because it is less sensitive to the learning rate than standard SGD.
The method is
detailed in \cref{alg:gcp-sgd}.
Based on the empirical results that follow, we recommend setting $s = d\bar n$.
If the total number of gradient samples is $s$, then we use
$s$ samples for uniform sampling in \cref{alg:stoc-grad-uniform} and
$p=\lfloor s/2 \rfloor$ and $q = \lceil s/2 \rceil$ for the
stratified sampling in \cref{alg:stoc-grad-stratified} and semi-stratified sampling
in \cref{alg:stoc-grad-semistrat}.
The parameter $\alpha$ is the learning rate and defaults to 0.01.
The Adam parameters are set to the values recommended in the original paper: $\beta_1 = 0.9$, $\beta_2 = 0.999$, and $\epsilon = 10^{-8}$.
We employ  a few standard modifications.
We group the iterates into epochs, and the number of iterations per epoch defaults to $\tau = 1000$.
To track progress, we estimate the function value $\FXM$ once per epoch.
Whenever the function value fails to decrease, we either decay the learning rate by $\nu = 0.1$ (as motivated by \cite{LoHu17}) or quit (after more than $\kappa=1$ failures).
We enforce any lower bounds $\ell$ given for the parameters by simple projection. In all of our examples, we use $\ell = 0$ (nonnegativity constraint).

\begin{algorithm}
  \caption{GCP with Adam (GCP-Adam)}
  \label{alg:gcp-sgd}\footnotesize
  \begin{algorithmic}[1]
    \Function{GcpAdam}{$\X$, $r$, $s$, $\alpha$, $\beta_1$, $\beta_2$, $\epsilon$, $\tau$, $\kappa$,  $\nu$,  $\ell$}
    \For{$k=1,2,\dots,d$}
    \State{\label{step:gcp-sgd-init}$\Ak \gets$ random matrix of size $n_k \times r$}
    \State $\Bk, \Ck \gets$ all-zero matrices of size $n_k \times r$
    \Comment{temporary variables used for Adam}
    \EndFor
    \State $\hat F \gets \textsc{EstObj}(\X, \Akset)$
    \Comment{estimate loss with fixed set of samples}
    \State $c \gets 0$
    \Comment{$c=$ \# of bad epochs (i.e., without improvement)}
    \State $t \gets 0$
    \Comment{$t=$ \# of Adam iterations}
    \While{$c \leq \kappa$}
    \Comment{$\kappa = $ max \# of bad epochs}
    \State Save copies of $\Akset$, $\Bkset$, $\Ckset$
    \Comment{save in case of failed epoch}
    \State $\hat F_{\text{old}} \gets \hat F$
    \Comment{save to check for failed epoch}
    \For{$\tau$ iterations}
    \Comment{$\tau = $ \# iterations per epoch}
    \State $\Gksset \gets \textsc{StocGrad}(\X, \Akset, s)$
    \Comment{$s = $ \# samples per stochastic gradient}
    \For{$k=1,\dots,d$}\tikzmark{AdamTop}
    \State $\Bk \gets \beta_1 \Bk + (1-\beta_1) \Gks$
    \State $\Ck \gets \beta_2 \Ck + (1-\beta_2) \Gks^2$
    \State $\Bk* \gets \Bk / (1-\beta_1^{t})$
    \State $\Ck* \gets \Ck / (1-\beta_2^{t})$
    \State $\Ak \gets \Ak - \alpha \, ( \Bk* \oslash \sqrt{\Ck*+\epsilon} )$\hspace{3mm}\tikzmark{AdamRight}\tikzmark{AdamBottom}
    \State $\Ak \gets \max\{\Ak, \ell\}$ \Comment{$\ell =$  lower bound}
    \EndFor
    \State $t \gets t+1$
    \EndFor
    \State $\hat F \gets \textsc{EstObj}(\X, \Akset)$  \Comment{estimate loss with fixed set of samples}
    \If{$\hat F > \hat F_{\text{old}}$} \Comment{check for failure to decrease loss}
    \State Restore saved copied of $\Akset$, $\Bkset$, $\Ckset$  \Comment{revert to last epoch's variables}
    \State $\hat F \gets \hat F_{\text{old}}$ \Comment{revert to prior function value}
    \State $t \gets t - \tau$ \Comment{wind back the iteration counter}
    \State $\alpha \gets \alpha \, \nu$ \Comment{reduce the learning rate}
    \State $c \gets c + 1$ \Comment{increment \# of bad epochs}
    \EndIf
    \EndWhile
    \State \Return $\Akset$
    \EndFunction
  \end{algorithmic}
  \AddNote{AdamTop}{AdamBottom}{AdamRight}{7cm}{Adam update depends on $\beta_1$, $\beta_2$, $\epsilon$; $\alpha =$ learning rate}
\end{algorithm}

We estimate $\FXM$ (via the function \textsc{EstObj}) in a way that is analogous to the stochastic gradient computation. There are two key differences. First, we use a much larger number of samples to ensure reasonable
accuracy, which is less important for the gradient computation. Second, the set of samples used for function estimation are fixed across all iterations for consistent evaluation across epochs.
For uniform sampling, let $\SIcnt$ be the number of times that index $i$ is selected and then estimate
\begin{equation}
  \label{eq:Fhat}
  \hat F \equiv \sum_{i \in \Omega}  \SIcnt \, \frac{n^d}{s} \, \fxme .
\end{equation}
For stratified sampling, let $\NIcnt$ denote the number of times that nonzero $i$ is selected and $\ZIcnt$ be the same for zero $i$ and then estimate
\begin{equation}
  \label{eq:Fhat-strat}
  \hat F \equiv \sum_{\xe \neq 0}  \NIcnt \, \frac{\eta}{p} \, \fxme
  + \sum_{\xe = 0}  \ZIcnt \, \frac{\zeta}{q} \, \fxme
  .
\end{equation}
In either case, it is easy to show that $\Exp{\hat F} = F$.
When we perform multiple runs of the same problem, we use the same set of samples for the function estimation \emph{across all runs} so that
they can be easily compared. The only exception is the non-stochastic method, which computes the full objective function.

\subsection{Sample Size and Comparison to Full Method for Dense Tensors}
\label{sec:sample-size}

We study the effect of sample size (also known as minibatch size) to understand the relevant trade-offs:
larger sample sizes yield lower variance stochastic gradients but higher costs per iteration.
We compute the GCP decomposition on an artificial four-way tensor of size $200 \times 150 \times 100\times 50$ and rank $r=5$
using the gamma loss function: $\fxm = x/m - \log m$ with a nonnegativity constraint on the factor matrices.
\Cref{sec:creat-gamma-test} provides the details of the data generation.

\begin{figure}
  \centering
  \subfloat[Individual runs with x-axis zoomed in to show the differences in the stochastic runs. \rundesc\@  \lossdesc]%
  {\label{fig:gamma-runs}\includegraphics[scale=0.5]{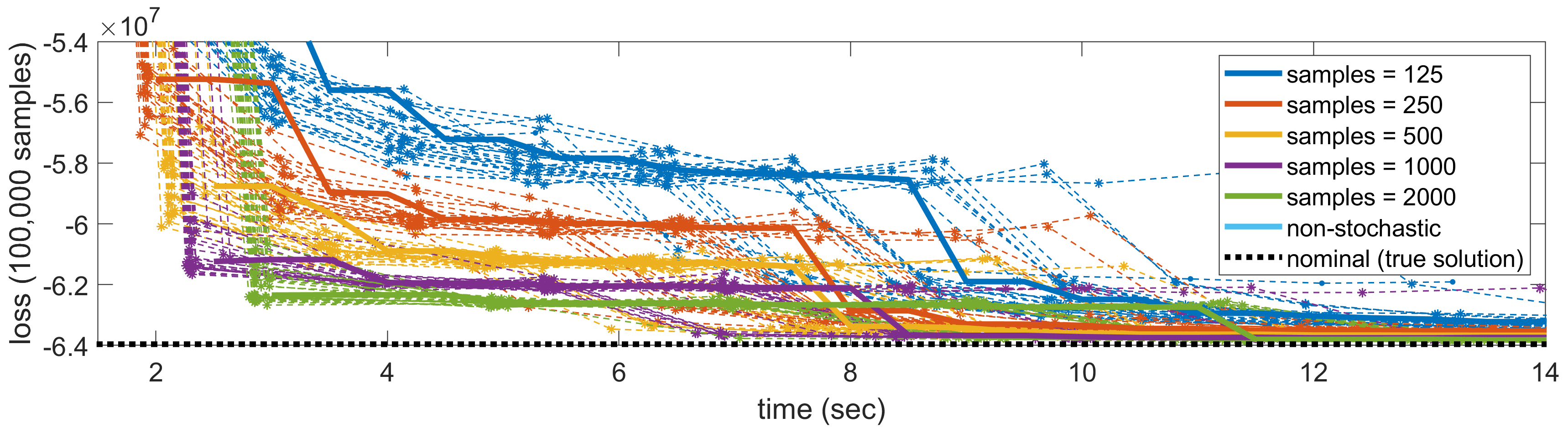}}
  \\
  \subfloat[Individual runs with x-axis zoomed out to show the
  non-stochastic method. For the non-stochastic method, each marker is
  a single iteration and the \emph{true} loss is plotted.]%
  {\label{fig:gamma-run-plus}\includegraphics[scale=0.5]{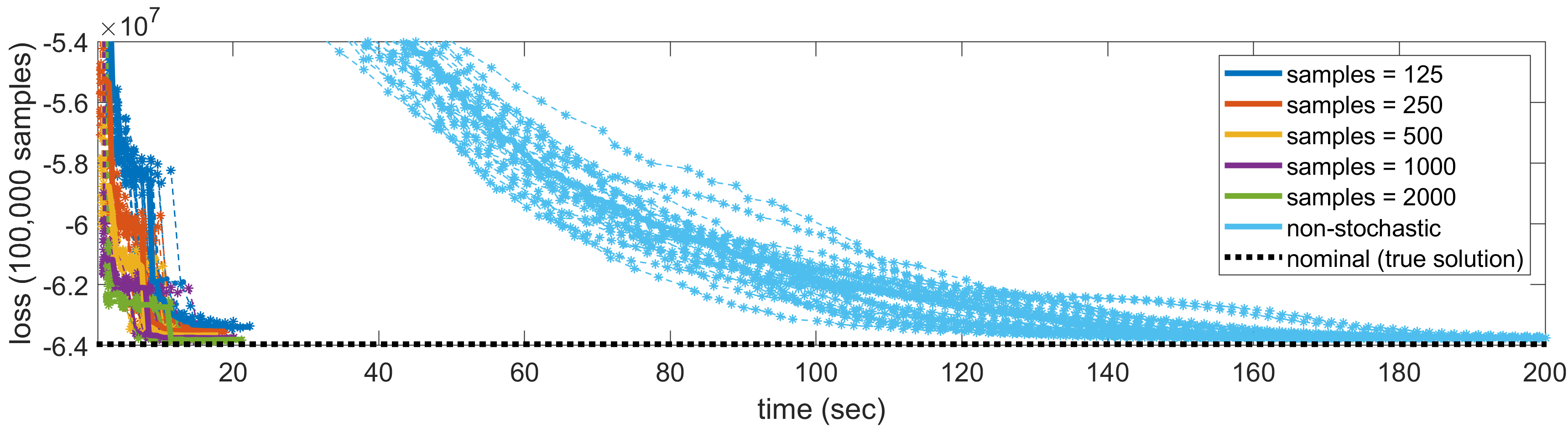}}
  \\
    \subfloat[Number of times the true solution was recovered, i.e., cosine similarity $\geq$ 0.9.]%
  {\label{fig:gamma-success}\includegraphics[scale=0.5]{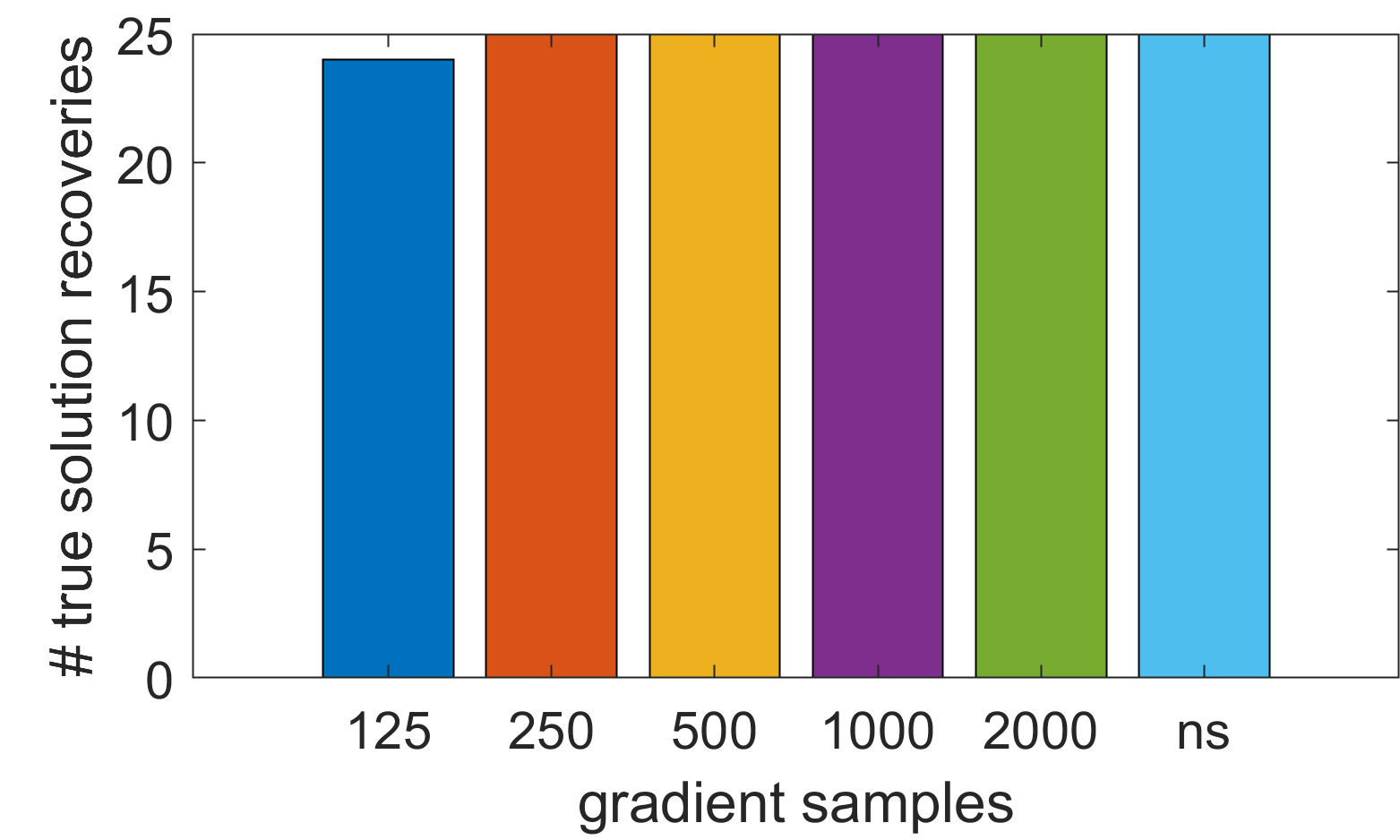}}
  ~~~~
  \subfloat[Boxplot of time per epoch. Each epoch is 1000 stochastic gradients plus one estimation of the function value.]%
  {\label{fig:gamma-epoch-time}\includegraphics[scale=0.5]{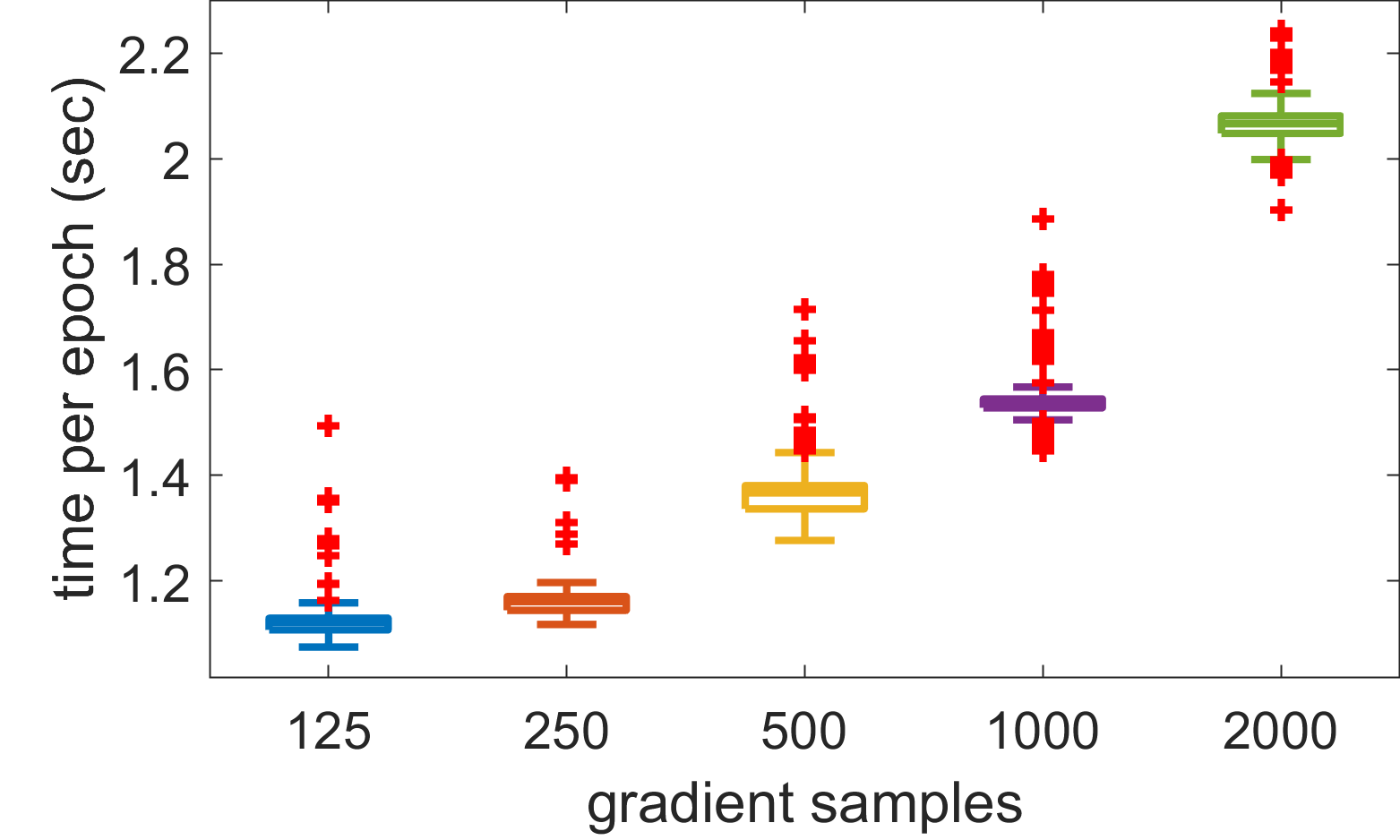}}
  \caption{GCP with Gamma loss $f(x,m) = x/m+\log m$ on artificial dense data tensor of size  $200 \times 150 \times  100\times  50$ and rank $r=5$, comparing various numbers of samples for the stochastic gradient in GCP-Adam and the non-stochastic GCP.
For each instance, we do 25 runs with different initial guesses. (The same 25 initial guesses are used for each instance.)}
  \label{fig:gamma}
\end{figure}

\Cref{fig:gamma} %
shows the results of GCP-Adam with various sample sizes ranging from $s=125$ to $s=2000$.
For comparison, we also include non-stochastic results based on the bound-constrained limited-memory BFGS method \cite{ByLuNoZh95}
using \emph{full} gradients; this optimization approach is standard in MATLAB toolboxes such as Tensor Toolbox \cite{AcDuKo11,TensorToolbox} and TensorLab \cite{VeDeDe16}.
The same set of 25 initial guesses is used for every instance. The initial guesses
comprise factor matrices with entries drawn uniformly from $(0,1)$.
For GCP-Adam, we estimate the loss $\hat F$ using 100,000 uniformly sampled entries that are fixed across all epochs and trials.
For the stochastic gradient, we use uniform sampling.

In the top two subfigures, we plot the function value versus time.
In \cref{fig:gamma-runs}, we consider just the stochastic methods and see the variation between them.
There are $d \bar n r = 2500$ free parameters and $\max_k n_k = 200$.
Common wisdom is to make one pass through the data per epoch, which
would require $s = n^d / \tau = 150,000$ samples per iteration.
However, we see that two orders of magnitude fewer samples are needed in practice, arguably due to the low-rank structure in the data.
Another option is to
set $s$ large enough so that we get at least one sample per row of $\Yks$ and thus every row of $\Ak$ is updated.
At a minimum, therefore, we may want $s \geq \max_k n_k = 200$.%
\footnote{This does not guarantee that every row is updated every time.
To do so with high probability, one might use roughly $10 \, \max_k n_k$ samples.
From \cite[Appendix A]{KoPiPlSe14}:
To sample $\rho t$  distinct ``types'' (i.e., row indices) from a set of $t$ types where $\rho \in (0,1)$, the
expected number of draws is $t \log(1/(1-\rho)) + O(1)$.
Therefore, the expected number of samples needed to collect 99.99\% of the members of a set with $t$ types is less than $10t$.}%
~%
Fewer samples generally corresponds to less progress per epoch;
however, sometimes fewer samples is advantageous because it
takes a different path to the solution.
For instance, we see that $s=2000$ initially makes better progress, but $s=500$ and $s=1000$ find the minimum more quickly.
At the other extreme, $s=125$ is the lowest cost per epoch, but its
progress in reducing the loss is hindered by too few samples.
In \cref{fig:gamma-run-plus} we show a longer time range on the x-axis so that the non-stochastic method is visible.
The non-stochastic method is approximately an order of magnitude slower.

The final function value may not tell the whole story, so we provide another measure of success.
Each instance is run with 25 random starts.
\Cref{fig:gamma-success} shows the number of times that each instance recovers the true solution,
meaning that the cosine similarity score between the true solution and the recovered solution is at least 0.9.
(See \cref{sec:meas-recov-true} for details of the cosine similarity score.)
The only time the true solution is \emph{not} recovered is one run for
$s=125$ samples per gradient.

\Cref{fig:gamma-epoch-time} shows box plots of
the cost per epoch, which is 1000 stochastic gradient evaluations (with $s$ specified on the x-axis) plus one function value estimation with 100,000 samples.
In each box plot, the middle line indicates the median, and the bottom and top edges of the box indicate the 25th and 75th percentiles, respectively. The whiskers extend to the most extreme data points not considered outliers, and the outliers are plotted individually as red '+' symbols.
Clearly, the cost per epoch grows linearly with the number of samples, but there are fixed costs that dominate the per iteration cost.
Notably, $s=2000$ uses $16$ times as many samples as $s=125$ but is generally only around twice as slow.

%

\subsection{Application to Dense Gas Measurements Tensor}
\label{sec:comparison-cp-als}

This section illustrates the benefits of using stochastic gradients
and the effect of sample size on GCP decomposition
for dense real data.
We consider a tensor based on chemo-sensing data collected
by Vergara et al.~\cite{VeFoMaTr13}.%
\footnote{Available at
  \url{http://archive.ics.uci.edu/ml/datasets/Gas+sensor+arrays+in+open+sampling+settings}.
  The dataset contains data for 11 gases; we used 7 that have more distinctive behaviors.
  It also has 6 sensor positions; we used position 3 (middle of the
  wind tunnel) where some of the interesting behaviors occur.
  The data is recorded at $\sim$100 Hz.
  We downsampled to 5 Hz,
  using the nearest measurement when needed,
  and skipped the first 9 seconds.
  We also removed sensor 33, which seemed to have erratic measurements.}
The dataset consists of measurements as a gas is blown over an array
of conductometric metal-oxide sensors
in a wind tunnel.
The tensor modes correspond to
 71 sensors,
 1250 time points,
 5 temperatures, and
 140 trials (7 gases with 20 repeats each).
This is a dense, relatively small 0.5~GB tensor.
We note that
Vervliet and De Lathauwer~\cite{VeLa16}
and Battaglino, Ballard and Kolda~\cite{BaBaKo18}
considered a similar 2~GB tensor derived from the same original dataset
with more time points but only three gases; however,
we do not compare directly to their approach because
they focus on standard
CP tensor decomposition.
For the loss, we use $\beta$-divergence with $\beta=1/2$ yielding the
loss function $f(x,m) = 2m^{1/2} + 2xm^{-1/2}$
with nonnegativity constraints.
This loss function is not necessarily optimal,
but it seemed to work well for this data in our experiments.
Furthermore, it is an attractive choice
since the tensor is nonnegative and may have some outliers.
We use rank $r=4$ because this is the smallest rank that
sufficiently distinguished the different gases in our experiments.

We run GCP-Adam  using uniform sampling and the non-stochastic GCP
under the same experimental conditions as in the previous subsection;
initial guesses are scaled to match the magnitude of the data tensor.
\begin{figure}
  \centering
  \includegraphics[scale=0.5]{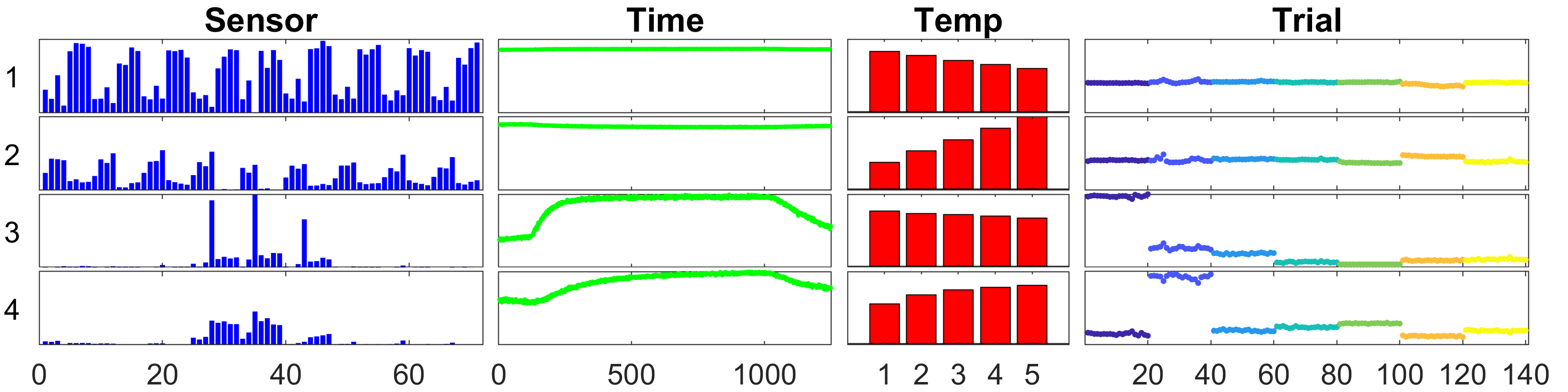}
  \caption{Factorization with the lowest overall score.
    GCP with $\beta$-divergence loss $f(x,m) = 2m^{\beta} + 2xm^{-\beta}$ with $\beta = 1/2$,
    rank $r=4$, and $s=1000$
    on
    real-world dense gas data tensor of size $71 \times 1250 \times 5 \times 140$.
    Sensor components are scaled proportional to component magnitude;
    the rest are normalized to length one.
    Trial symbols are color coded by gas.
  }
  \label{fig:gas-viz}
\end{figure}
\Cref{fig:gas-viz} shows the results of the best overall run in terms of the
final objective value.
The components are ordered by magnitude,
the sensor mode is normalized to the magnitude of the component,
and the other modes are normalized to unit length.
The trial mode is color coded by gas.
The first two components focus largely on
sensor variations due to temperature.
The final two components capture
some temporal patterns for each gas
that impact sensors near the middle
of the array,
where the gas is likely most concentrated.
The factorization identifies generally smooth temporal profiles
and tends to group the same gas (indicated by color) in the trial mode.

\begin{figure}
  \centering
    \subfloat[Boxplots of total runtimes.
    \label{fig:gas-sample-size-runtimes}]
    {\includegraphics[scale=0.5]{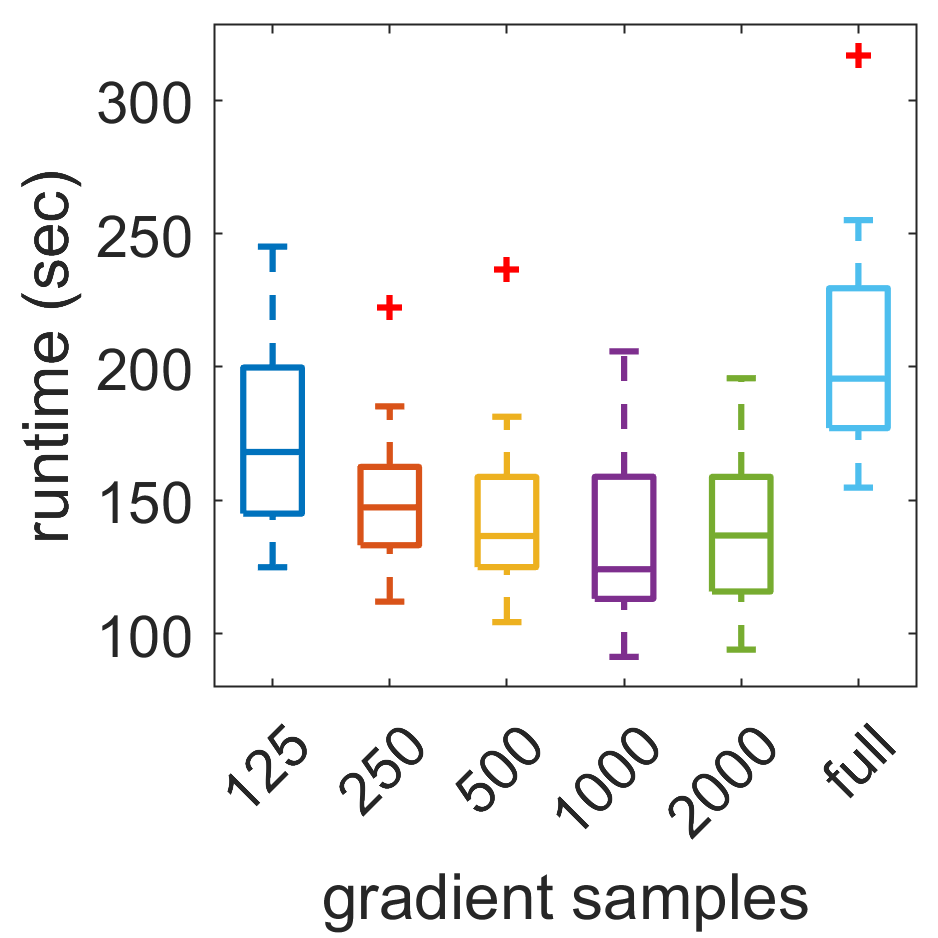}}
    ~~~~
  \subfloat[Boxplots of final model losses.
    \label{fig:gas-sample-size-loss}]
    {\includegraphics[scale=0.5]{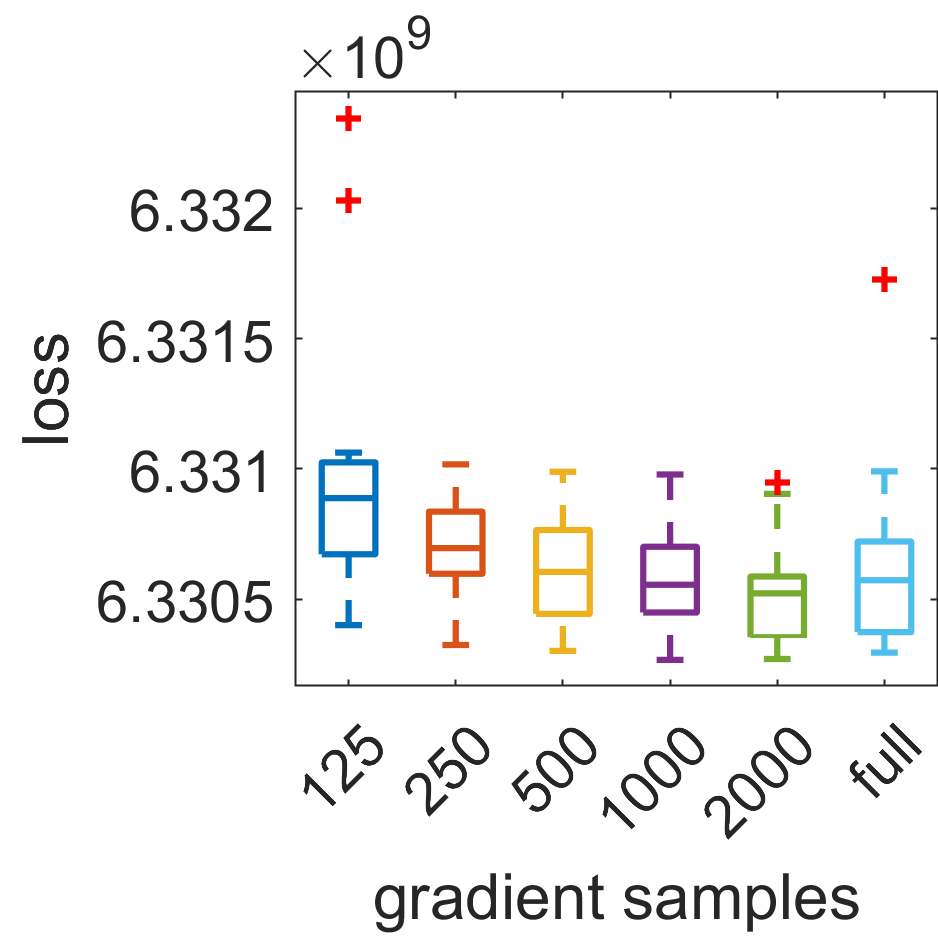}}
    ~~~~
    \subfloat[Clustering performance of each run. We apply k-means to the fourth factor matrix and measure what percentage of the 140 trials are clustered correctly by gas.]
    {\label{fig:gas-sample-size-clustering}%
      \includegraphics[scale=0.5]{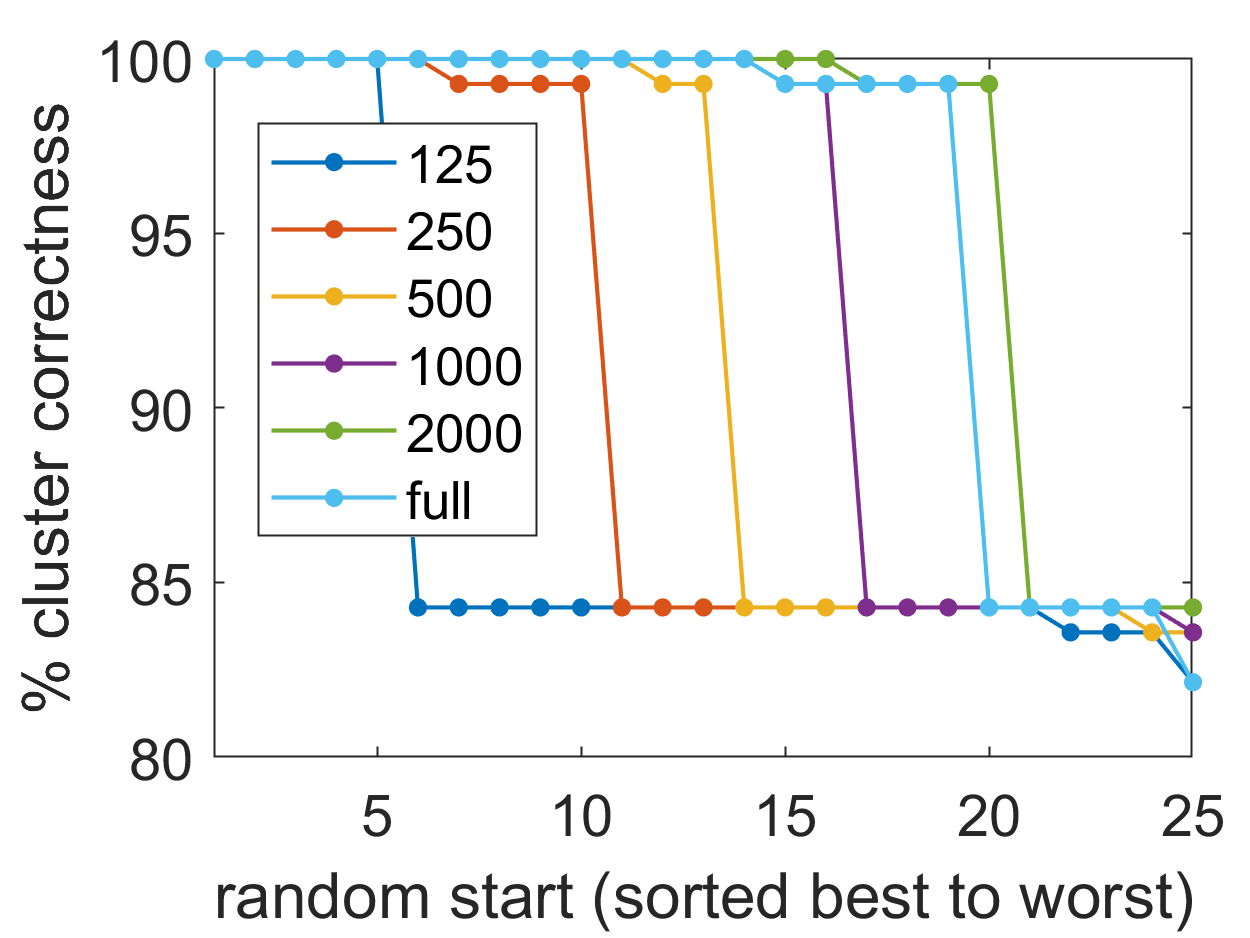}}
    \caption{Comparison of
      GCP-Adam with uniform sampling and sample sizes from $s=125$ to $s=2000$ as well as
      non-stochastic GCP (``full''),
      for fitting a dense gas tensor of size $71 \times 1250 \times 5 \times 140$
    with $\beta$-divergence loss using $\beta = 1/2$, i.e., $f(x,m) = 2m^{1/2} + 2xm^{-1/2}$.
    Each box plot is based on 25 runs with different initial guesses.
    (The same 25 initial guesses are used for each box plot.)}
  \label{fig:gas-various}
\end{figure}

\Cref{fig:gas-various} shows performance results.
\Cref{fig:gas-sample-size-runtimes} shows boxplots of the runtimes
for the varying sample sizes and the non-stochastic GCP (``full'').
The smaller sample sizes take longer because they converge more slowly even though
each epoch is cheaper.
The non-stochastic instance is overall slowest.
\Cref{fig:gas-sample-size-loss} shows the range of objective function values
for each instance, which is very small.
For the stochastic instances, larger sample sizes tend to achieve a marginally better loss,
as indicated by improved median and percentile losses.
The stochastic instances with $s \geq 500$ samples per gradient perform at least as well as the non-stochastic instance,
while being much faster.
\Cref{fig:gas-sample-size-clustering}
shows performance on a clustering task.
For each run (i.e., a given random start and instance),
the rows of the fourth factor matrix are clustered via k-means%
\footnote{We used 500 replicates with MATLAB's built-in \texttt{kmeans} command
  to avoid local minima.}
and we measure what percentage of the 140 trials (7 gases with 20 trials each)
get clustered correctly.
Larger sample sizes have generally better clustering performance here.
Non-stochastic GCP performs
similarly to $s = 2000$ samples per gradient.

%
%
%
%

\subsection{Sample Size for Sparse Tensors}
\label{sec:sample-size-sparse}

We consider a four-way \emph{sparse} binary tensor of size  $200 \times 150 \times 100\times 50$ and rank $r=5$ generated
according to an \emph{odds} model tensor, i.e., where $\me$ gives the odds that $\xe=1$.
The procedure to generate the data is described in detail in \cref{sec:creating-binary-test}.
The factor matrices in the solution $\M$ have $(r-1)$ sparse columns and one constant column.
The result is a tensor that has 150,452 `structural' nonzeros (from the sparse columns) and 374,435 `noise' nonzeros (from the dense column),
with an overall density of 0.35\%.
We use the loss corresponding to Bernoulli data with an odds link, i.e., $f(x,m) = \log(m+1) - x \log m$
and a nonnegativity constraint on the factor matrices.

\begin{figure}
  \centering
  \subfloat[Individual runs. \rundesc\@  \lossdesc]%
  {\label{fig:binary-runs}\includegraphics[scale=0.5]{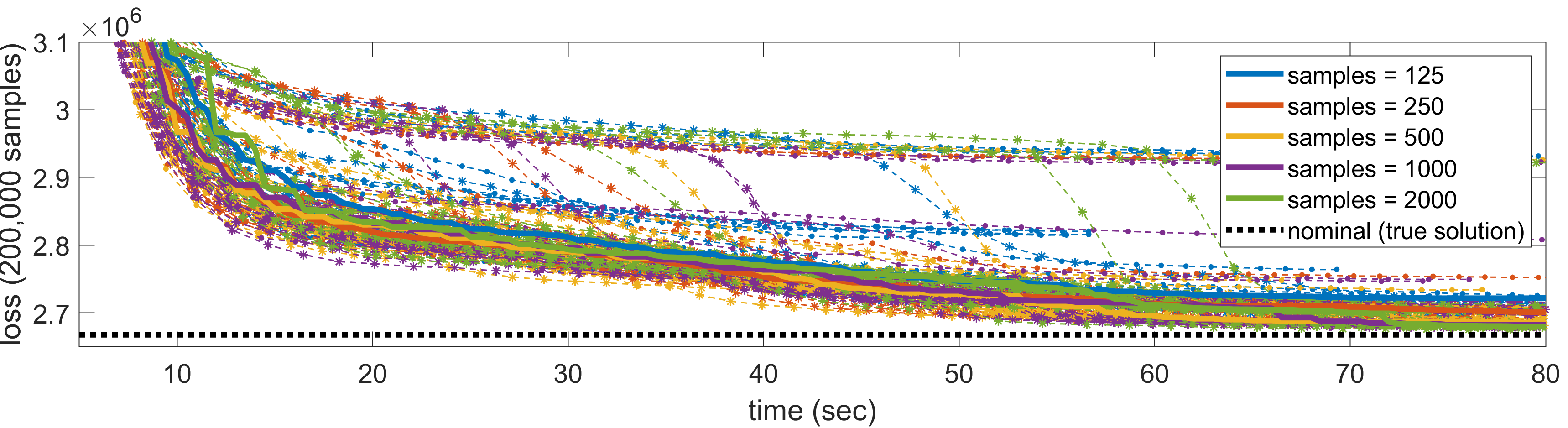}}
  \\
  \subfloat[Number of times the true solution was recovered, i.e., cosine similarity $\geq$ 0.9.]%
  {\label{fig:binary-success}\includegraphics[scale=0.5]{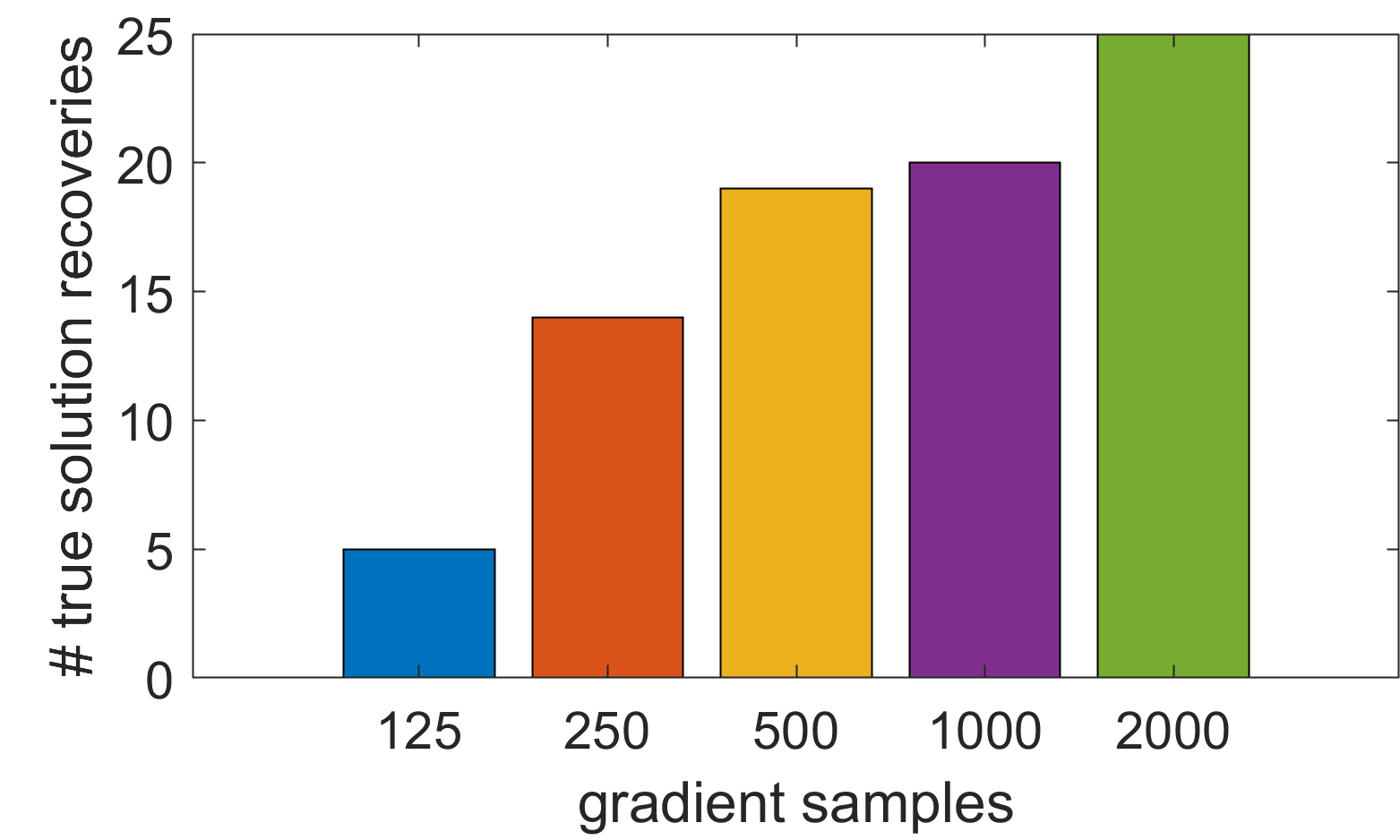}}
  ~~~~
  \subfloat[Boxplot of time per epoch. Each epoch is 1000 stochastic gradients plus one estimation of the function value.]%
  {\label{fig:binary-epoch-time}\includegraphics[scale=0.5]{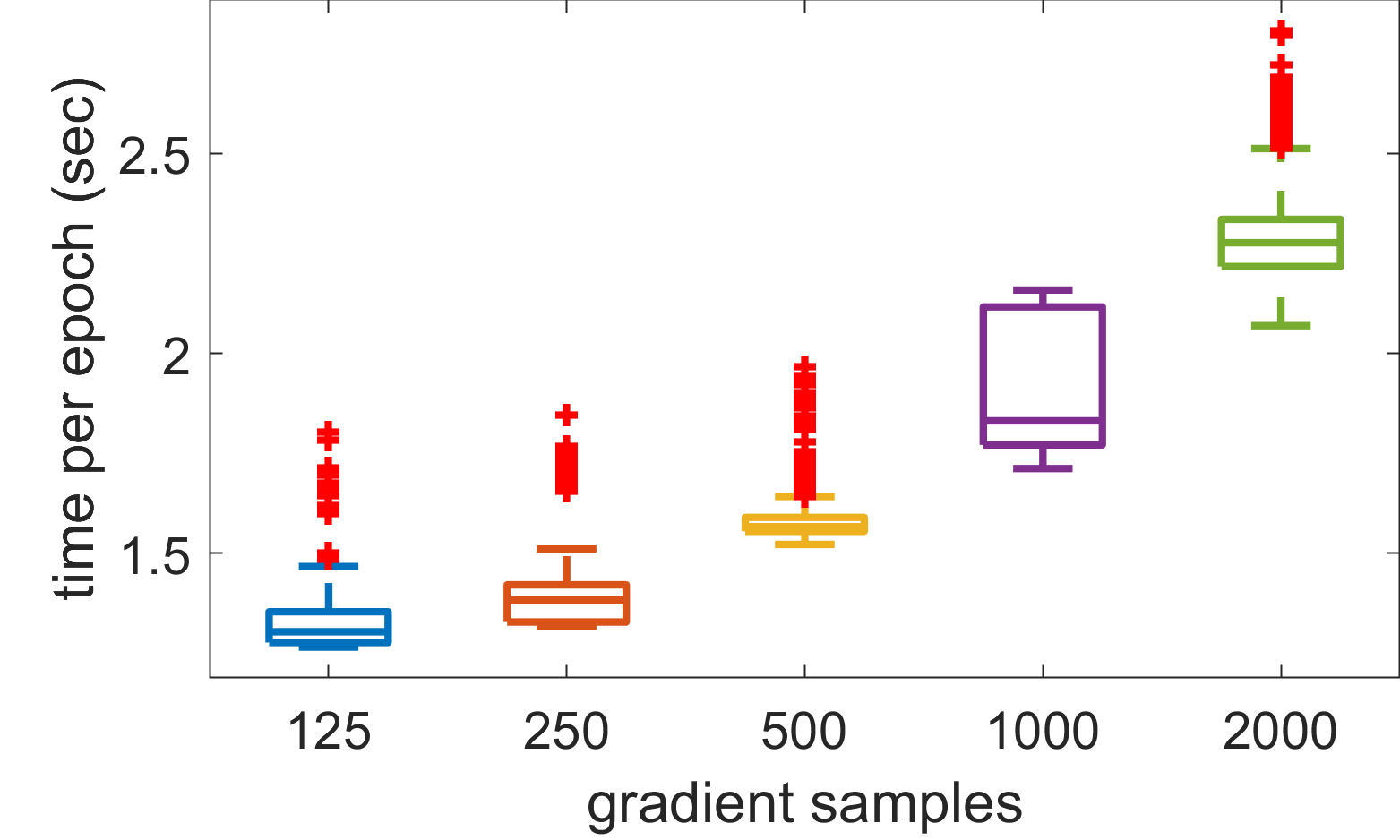}}
  \caption{GCP with Bernoulli loss $f(x,m) = \log(m+1) - x \log m$ on artificial sparse data tensor of size  $200 \times 150 \times  100\times  50$ and rank $r=5$ with 524,468 (0.35\%) nonzeros. Comparing various numbers of samples for the stochastic gradient in GCP-Adam with stratified sampling.
For each instance, we do 25 runs with different initial guesses. (The same 25 initial guesses are used for each instance.)}
  \label{fig:binary}
\end{figure}

The results of GCP-Adam with various sample sizes ranging from $s=125$ to $s=2000$ are shown in \cref{fig:binary}.
We use stratified sampling to compute the stochastic gradient,
and the gradient samples are evenly divided between zeros and nonzeros.
The same set of 25 initial guesses is used for every instance. The initial guesses
comprise factor matrices with entries drawn uniformly from $(0,1)$ and then scaled to match the magnitude of the true solution tensor.
For GCP-Adam, we estimate the loss $\hat F$ using 200,000 stratified sampled entries (evenly divided between zeros and nonzeros) that are fixed across all epochs and trials.

\Cref{fig:binary-runs} plots the individual runs.
For the stochastic method, the overall run time is about four times more than for the (easier) Gamma case.
We omit the non-stochastic method since it is again significantly slower.
This is arguably a difficult test problem in terms of recovering the true factors, especially compared to the Gamma problem in
\cref{sec:sample-size}.
From \cref{fig:binary-success}, observe that $s=125$ fails to find the true solution more often than it succeeds.
For $s \geq 250$, the true solution is recovered in the majority of cases, and the recovery rate improves as the number of samples increases;
all $25$ runs in this experiment succeeded for $s=2000$.
\Cref{fig:binary-epoch-time} shows boxplots of the time per epoch.
As was the case with uniform sampling of a same-sized dense tensor in \cref{sec:sample-size},
the time per epoch mainly consists of costs that grow linearly with the number of samples
and fixed costs that, at this scale, remain significant.
As before, $s=2000$ uses $16$ times as many samples as $s=125$
but is generally only around twice as slow.

\subsection{Comparison of Sampling Strategies for Sparse Tensors}
\label{sec:comp-sampl-strat}

This section shows the benefit of using stratified and semi-stratified sampling over uniform.
We use the same procedure as \cref{sec:sample-size-sparse} to create a sparse binary tensor (detailed in \cref{sec:creating-binary-test}).
We generate a tensor of size $400 \times 300 \times 200 \times 100$ that is 0.38\% dense.
This example has 4,402,374 `structural' nonzeros and 4,788,052 `noise' nonzeros, with a total of 9,181,549 nonzeros (less than the sum due to overlap).
Storing it takes 0.37~GB as a sparse tensor,
but would take 19~GB as a dense tensor.

\begin{figure}
  \centering
  \subfloat[Individual runs. Each dashed line represents a single run, and the markers signify epochs. The marker is an asterisk  if the true solution was recovered and a dot otherwise. Solid lines represent the median. The dashed black line is the function value estimate for the true solution. Across all runs, the \emph{same} set of samples is used to estimate the loss.]%
  {\includegraphics[scale=0.5]{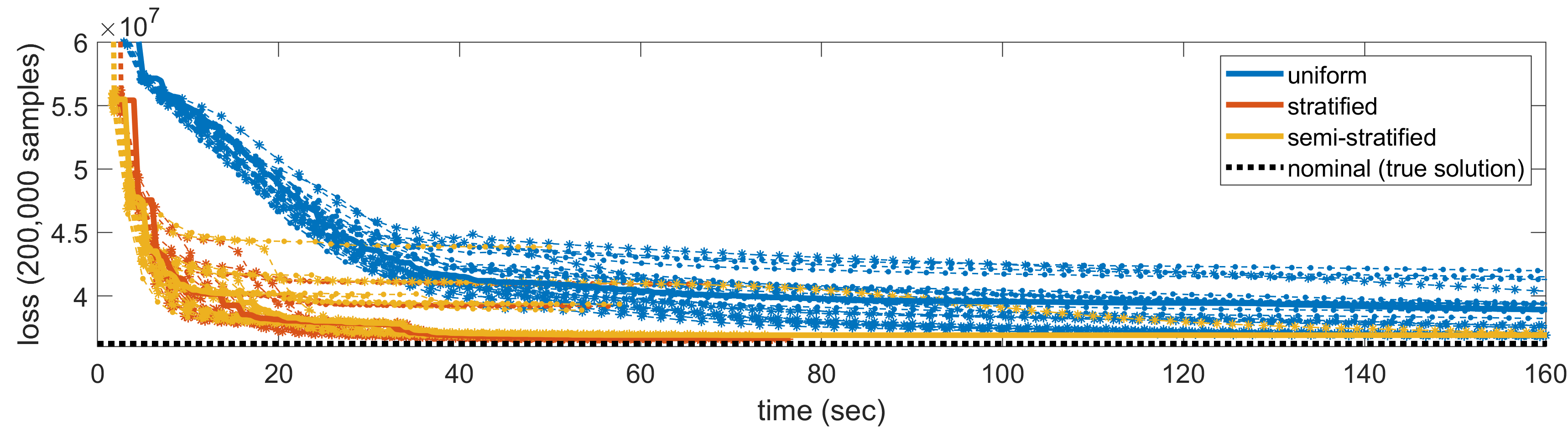}}\\
 \subfloat[Number of times the true solution was recovered, i.e., cosine similarity $\geq$ 0.9, for each number of gradient samples.]{
   ~~~~\includegraphics[scale=0.5]{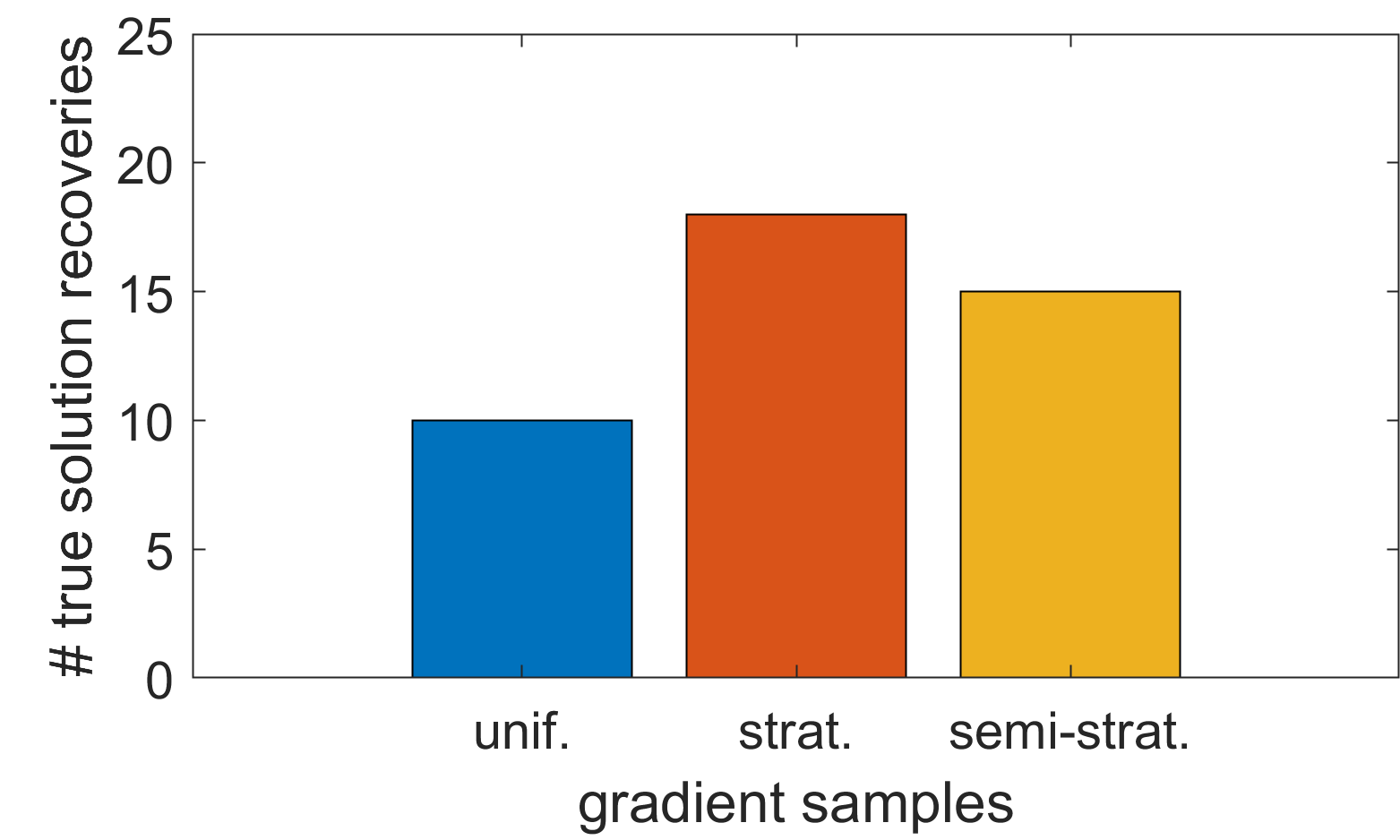}~~}~~
 \subfloat[Box plot of time per epoch for the different sampling methods, with the midline representing the median time.]{
  ~~\includegraphics[scale=0.5]{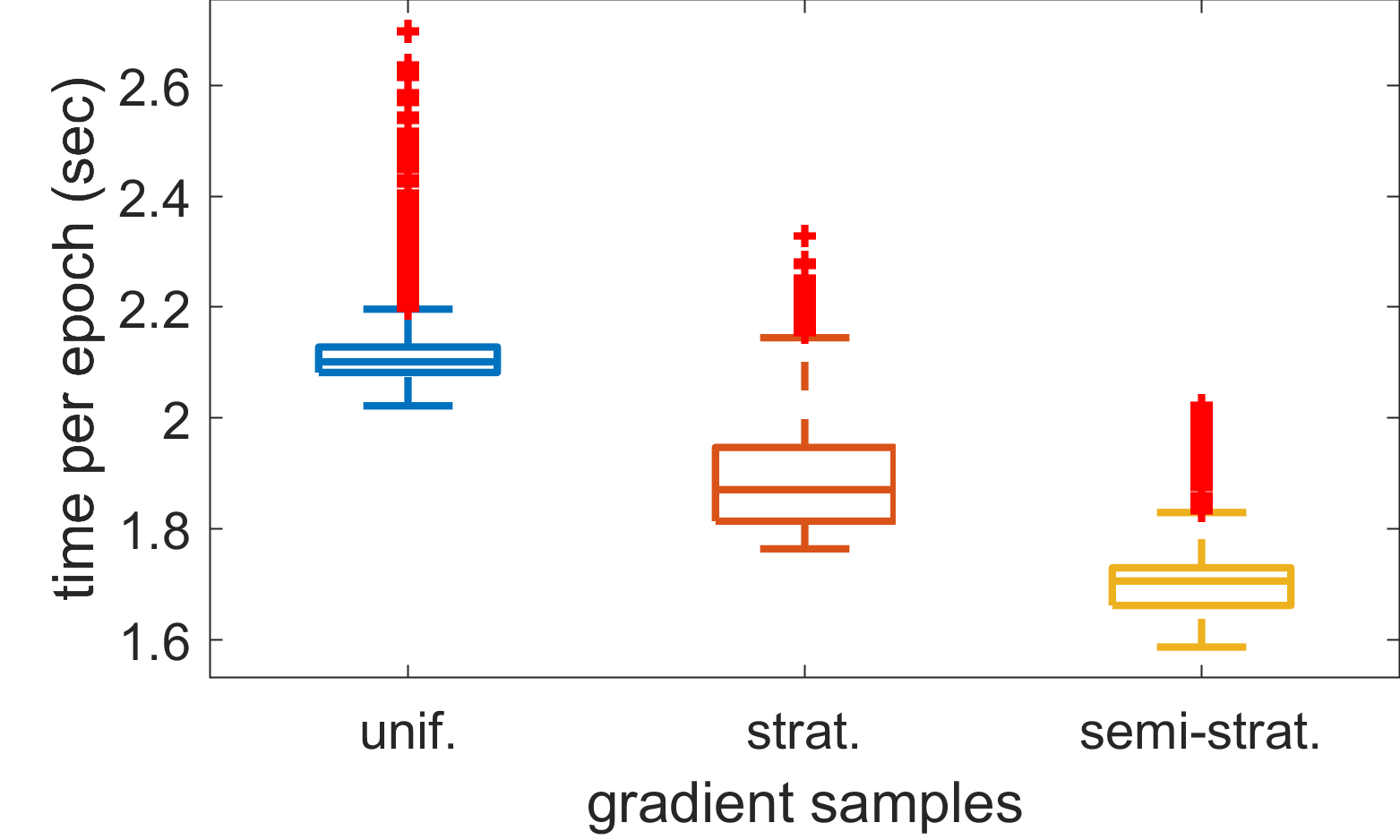}~~~~}
\caption{GCP with Bernoulli loss $f(x,m) = \log(m+1) - x \log m$ on artificial sparse data tensor of size $400 \times 300 \times 200 \times 100$ with rank $r=5$ and 9,181,549 nonzeros (0.38\% dense). Comparing different sampling strategies in GCP-Adam with $s=1000$ samples per stochastic gradient, evenly divided between nonzeros and zeros (or ``zeros'') for the stratified samplers. For each strategy, we do 25 runs with different initial guesses. (The same 25 initial guesses are used for each instance.)}
  \label{fig:compare}
\end{figure}

The results of uniform, stratified, and semi-stratified sampling are shown in \cref{fig:compare}.
We calculate the estimated loss $\hat F$ once per epoch using 200,000 stratified sampled entries, evenly divided between zeros and nonzeros.
We use $s=1000$ samples per stochastic gradient evaluation ($s = d\bar n$), evenly divided between nonzeros and zeros (or ``zeros'') for the stratified samplers.
We use 25 initial guesses with random positive values, scaled so that the norm of the initial guess is the same as the norm of the tensor.
Note that uniform sampling on a sparse tensor is not as fast as it is on a dense tensor since only nonzeros are stored explicitly and
every sampled index has to be checked against the list of nonzeros to determine its value.

The proportion of nonzeros is arguably higher than what is observed in many real-world datasets,
which is a favorable condition for uniform sampling since every sample will more likely include nonzeros.
Nevertheless,
stratified and semi-stratified approaches clearly outperform uniform sampling.
They converge faster and more often find the true solution, all at less cost per epoch.
This result is expected because stratified sampling should reduce the variance.
The time advantage of the semi-stratified approach is minimal for this small example.
However, the speed of the MATLAB implementation of stratified sampling depends on the ability to use linearized indices, which means that the total size of the tensor must be less than $2^{64}$. For larger tensors where we cannot use this approach, the sampling efficiency can degrade by more than an order of magnitude.
The semi-stratified approach has no such limitation.

\Cref{tbl:sampling:variance} provides empirical validation that
the variances of the stratified and semi-stratified sampling are lower than for uniform.
We consider the variances with respect to the  vectorized gradient, i.e.,
$\gvec = [ \vc(\Gk{1});\, \vc(\Gk{2});\, \dots;\, \vc(\Gk{d}) ] \in \Real^{d \bar n r}$.
Let $\gvecs$ denote the random variable representing the vectorized stochastic gradient, and recall that $\Exp{\gvecs} = \gvec$
because the sampling methods lead to unbiased estimators.
Letting $\gvecs_{\xi}$ denote the $\xi$-th of $N$ realizations of the random variable, the quantities of interest are
\begin{displaymath}
  \text{empirical bias } = \| \V[\hat]{g} - \gvec \|_2
  \qtext{where} \V[\hat]{g} = \frac{1}{N} \sum_{\xi=1}^N \gvecs_{\xi},
  \quad\text{and}
\end{displaymath}
\begin{displaymath}
  \text{empirical variance } =
  \text{trace}\left(
    \frac{1}{N} \sum_{\xi=1}^N (\gvecs_{\xi}- \V[\hat]{g})
    (\gvecs_{\xi}- \V[\hat]{g})^{\intercal}
  \right)
  = \frac{1}{N} \sum_{\xi=1}^N \| \gvecs_{\xi}- \V[\hat]{g}\|_2^2.
\end{displaymath}
\Cref{tbl:sampling:variance} considers two distinct cases:
an initial guess (far from the solution with a larger gradient)
and the overall best final solution (close to the solution with a smaller gradient).
In both cases,
stratified and semi-stratified sampling had lower variances than uniform.
Consequently, these have lower empirical biases as well.
This helps to explain the superior convergence of the stratified and
semi-stratified approaches.

\begin{table}
  \centering\footnotesize
  \caption{Empirical bias and variance of stochastic gradients at an initial guess and a final model for GCP with Bernoulli loss using the artificial sparse tensor and sampling methods of \cref{fig:compare}. For each method and model, we use $N=1000$ stochastic gradient realizations.}
  \begin{tabular}{|l|rc|rc|}
  \hline
sampling               & \multicolumn{2}{c|}{Initial guess: $\|\gvec\|_2 = $ 2.00e+07}  & \multicolumn{2}{c|}{Final model: $\|\gvec\|_2 = $ 3.10e+06} \\

method                  & ~~~~emp.~bias                          & emp.~var.                      & ~~~~emp.~bias                          & emp.~var.       \\
  \hline
  uniform         & 1.08e+06                & 1.26e+15                & 3.90e+07                & 1.52e+18 \\
  stratified      & 7.48e+05                & 5.64e+14                & 3.17e+06                & 9.84e+15 \\
  semi-stratified & 7.68e+05                & 5.70e+14                & 3.14e+06                & 9.93e+15 \\
  \hline
  \end{tabular}
  \label{tbl:sampling:variance}
\end{table}

\subsection{Application to Sparse Count Crime Data and Comparison to CP-APR}
\label{sec:application-crime}

We consider a real-world crime statistics dataset comprising more than 15 years of crime data from the city of Chicago.
The data is available at \url{www.cityofchicago.org}, and we downloaded a 4-way tensor version from
FROSTT~\cite{frosttdataset}. The tensor modes correspond to
 6,186 days from 2001 to 2017,
 24 hours per day,
 77 communities, and
 32 types of crimes.
Each $\X(i,j,k,\ell)$ is the number of times that crime $\ell$ happened in neighborhood $k$ during hour $j$ on day $i$.
The tensor has 5,330,673 nonzeros. Stored as a sparse tensor, it requires 0.21~GB of storage.

Since this is count data, we use GCP with a Poisson loss function, i.e., $f(x,m) = m - x \log m$ and nonnegativity constraints.
We run GCP-Adam with both stratified and semi-stratified sampling using $s=d \bar n = 6319$ samples for the stochastic gradient.
We compare to  the state-of-the-art for CP alternating Poisson Regression (CP-APR) \cite{ChKo12},
using the Quasi-Newton method described in \cite{HaPlKo15}. We run each method with 20 different starting points.
We compute rank $r=10$ factorizations.

\begin{figure}
  \centering
  \includegraphics[scale=0.5]{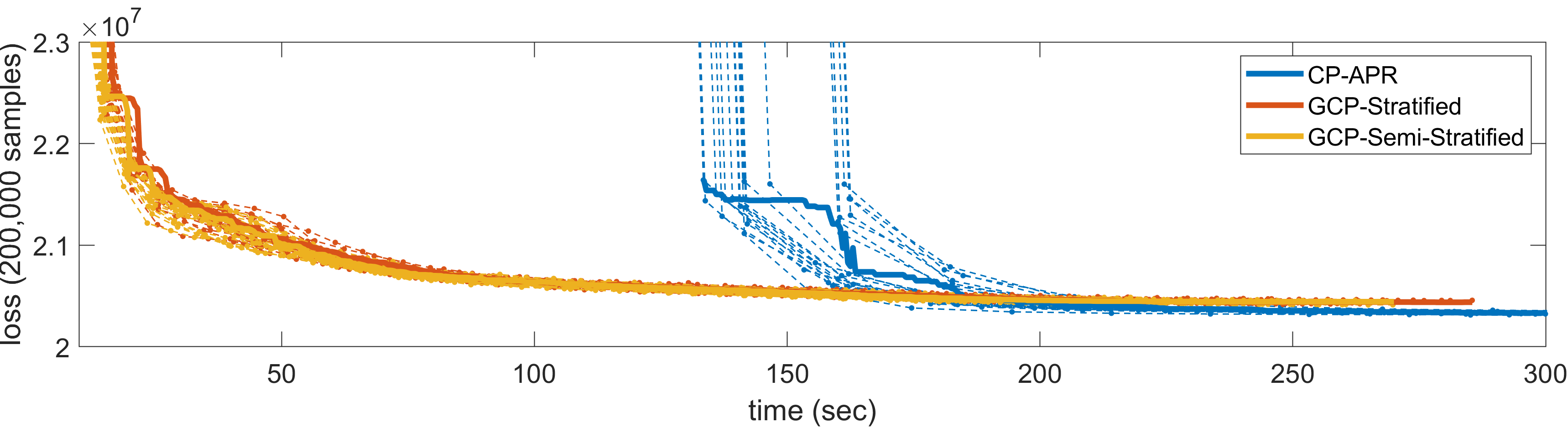}
  \caption{GCP with Poisson loss $f(x,m) = m - x \log m$, rank $r=10$, and $s=6319$ (sum of the dimensions) on real-world Chicago crime
    data tensor of size $6186 \times 24 \times 77 \times 32$ with
    5,330,673 nonzeros. Comparing GCP-Adam using both stratified and semi-stratified sampling
    with CP-APR (using quasi-Newton).
    We run each method from 20 different initial guesses.
    Each dashed line represents a single run, and the markers signify epochs for GCP and iterations for CP-APR.
    Solid lines represent the median.  Across all GCP runs, the \emph{same} set of samples is used to estimate the loss.
  CP-APR computes the exact loss, and the differences between the final losses seem to be in part an artifact of the estimation.}
  \label{fig:chicago-compare}
\end{figure}

Timing results are shown in \cref{fig:chicago-compare}. The CP-APR method has to do extensive pre-processing, which is why the
first iteration does not complete until approximately 140 seconds.
The GCP-Adam methods descend much more quickly but do not reduce the loss quite as much,
though this failure to achieve the same final minima is likely an artifact of the function estimation and/or nuance of the Adam parameters.

\begin{figure}
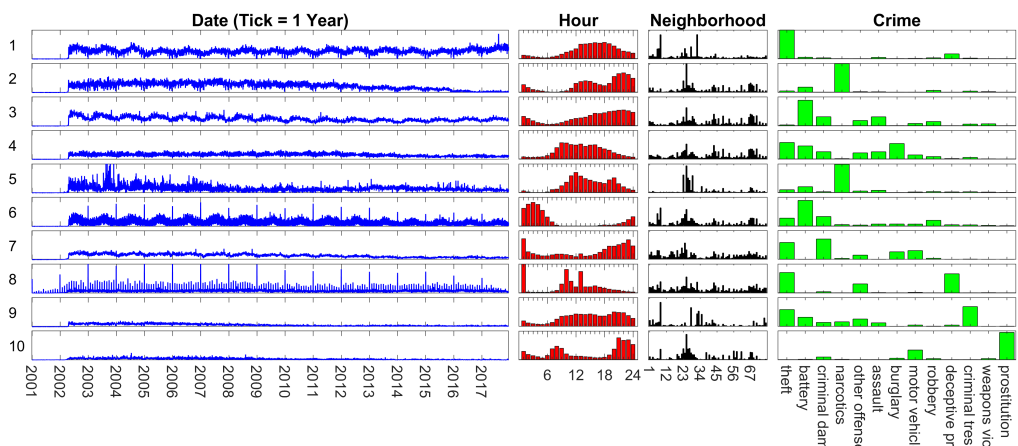
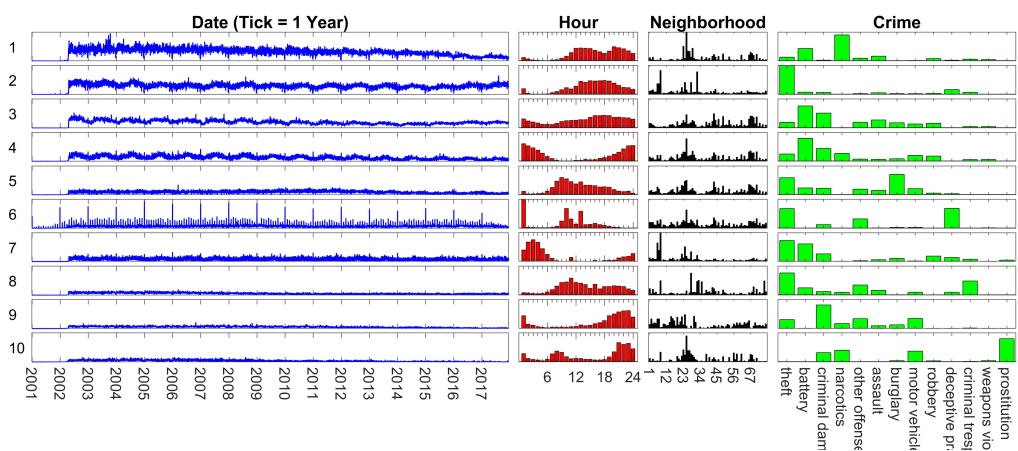
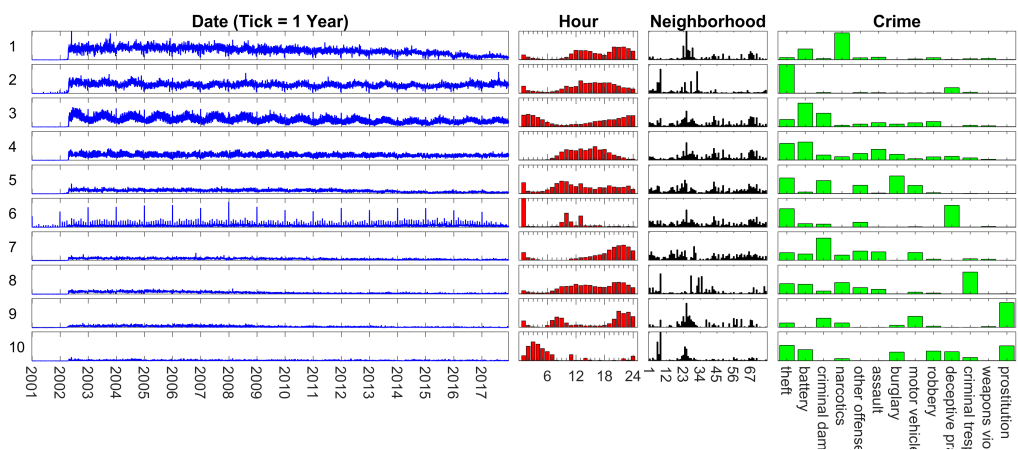

  \centering
  \subfloat[CP-APR (non-stochastic)]{\cpic{cpapr}}\\
  \subfloat[GCP-Adam with stratified sampler]{\cpic{strat}}\\
  \subfloat[GCP-Adam with semi-stratified sampler]{\label{fig:chicago-semistrat}\cpic{semistrat}}
  \caption{Visualization of factorizations of the Chicago crime tensor as produced by three different methods.}
  \label{fig:chicago-viz}
\end{figure}

Although the final loss  functions values are slightly different,
the factorizations computed by the three methods are similar.
We show the results from the first random starting point in \cref{fig:chicago-viz}.
Each rank-1 component corresponds to a row in the figure.
The first column is the day, shown as a line graph.
Components in this column are scaled to capture the magnitude of the overall component,
while all the other components are normalized.
The second column is hour of the day, shown as a bar graph.
The third column is the neighborhood, shown as a bar graph.
(The third column is somewhat difficult to interpret visualized this way but we show it on a map in \cref{sec:indiv-comp-chic}.)
The fourth column is the crime type, sorted by overall frequency, and only showing the 13 most prevalent crimes.
When working with real-world data, we generally have to experiment.
Nevertheless, certain trends emerge over and over again for different starting points, ranks, and methods.
In this case, we see strong commonalities among the three methods, as follows.
Component 1 for CP-APR is similar to component 2 for the GCP-Adam methods, with ``theft'' being the most prevalent crime and a similar pattern in time.
Component 2 for CP-APR is similar to component 1 for the GCP-Adam methods.
Component 3 for the semi-stratified GCP-Adam has a very strong seasonal signature, becoming most active in summer months.
  Component 4 for stratified GCP-Adam is similar, as is component 6 for CP-APR.
 Component 8 in CP-APR and component 6 in both GCP-Adam methods has a special pattern of a spike on the first of each year and again on the first of each month. There is also a spike at midnight in the hour mode. This is likely a feature of how the associated crimes were recorded in the dataset.
To help further with interpretation of the rank-1 components, we zoom in on the components from the semi-stratified GCP-Adam solution
in \cref{sec:indiv-comp-chic}, where we show each individual component, including drawing a heatmap of the neighborhoods on a map.

\section{Conclusions}
\label{sec:conclusions}

We propose a stochastic gradient for GCP tensor decomposition
with general loss functions.
The structure of the GCP gradient means that there is no general way
to maintain sparsity in its computation even when the input
tensor is sparse.
Our investigation was prompted by the sparse case,
but the stochastic approach applies equally well to dense tensors.
A unique feature of our approach is the use of stratified
and semi-stratified
sampling in the gradient computation for sparse tensors.

We tested the stochastic gradient using Adam \cite{KiBa15} for GCP tensor decomposition
and made several findings.
Stochastic gradient methods are effective in practice in terms of driving down the objective function and recovering the true solutions.
Empirically, we find that the number of samples should be roughly equal to the sum of the dimensions, i.e., $s = \sum_{k=1}^d n_k = d \bar n$.
This is much lower than would be required to cycle through the entire dataset each epoch,
i.e., $n^d$ divided by the number of iterations per epoch.
For dense problems, stochastic gradient methods can be much faster than the non-stochastic prior approach using L-BFGS-B \cite{HoKoDu20}.
For sparse problems, stochastic gradients enable us to circumvent formation of the dense tensor needed by the gradient, making it possible to solve much larger problems.
Additionally, we propose stratified and semi-stratified sampling, which are typically superior to uniform sampling.
We have not discussed how to determine the rank since that is a difficult problem even
for standard CP.

Overall, stochastic methods have proved to be a promising approach for GCP tensor decomposition,
especially for large-scale sparse tensors which have no viable alternative.
However, many open questions remain.
We can likely further improve the results by using more sophisticated stochastic optimization methods,
e.g., weight decay strategies in Adam \cite{LoHu17}.
Likewise, more sophisticated sampling strategies such as leverage score sampling from matrix sketching \cite{Ma11,Wo14,ChBhSaWa15}
may further improve performance by reducing the variance.
Another important line of investigation is developing appropriate theory to describe the
improvement gains of the stratified approaches.

In terms of implementations, an interesting consequence of sampling in the context of parallel tensor decomposition~\cite{SmRaSiKa15,KaUc15,LiChPeSu17,PhKo19}
is that we can reduce the computation and/or communication by sampling only a subset of the entries.
Moreover, we may be able to stratify the samples in such a way that is amenable to more structured communications.

\appendix

\section{Special Cases where Gradient Does Not Require Dense Calculations}%
\label{sec:special-cases}%
Computing the gradient is not a major issue for standard CP due to its special structure.
Specifically, the computation of $\Gk$ can be simplified so that
the primary work is computing $\Xk\Zk$, which is a sparse MTTKRP whenever $\X$ is sparse \cite[Appendix A]{HoKoDu20}.
For Poisson CP~\cite{ChKo12}, the primary work is computing $\Mx{V}_{\!(k)}\Zk$ where $\T{V}$ is the sparse tensor defined as
\begin{displaymath}
  v_i =
  \begin{cases}
    \xe / \me & \text{if } \xe \neq 0, \\
    0 & \text{otherwise}.
  \end{cases}
\end{displaymath}
In these cases, the gradient can be computed in $O(\nzX\,rd)$ flops with  $O(\nzX)$  additional storage.
Generally, however, we may not have such structure and we have to compute with a dense $\Y$ tensor
at a cost of $O(rn^d)$ flops and $O(\nzX)$ extra storage.

\section{Determining the Oversampling Rate}%
\label{sec:determ-overs-rate}%
\Cref{sec:stratified-sampling-1} mentions that we oversample to get sufficiently many zeros with high probability.
Namely, we sample
\begin{equation*}
\rho \, \frac{n^d}{n^d-\nzX} \, \szero
= \rho \, \frac{1}{1-\nzX/n^d} \, \szero
= \rho \, \frac{\szero}{\pzero}
\end{equation*}
indices to get the desired $\szero$ zeros,
where $\pzero = 1 - \nzX/n^d$ is the proportion of zeros in the tensor.
Here we discuss the oversampling rate $\rho$.

We can use the inverse cumulative distribution function (CDF) of a negative binomial distribution to determine an appropriate $\rho$.
The negative binomial distribution models the number of failures before a given number of successes with a specified success rate.
In our case, we want $\szero$ successes %
and the success rate is $\pzero$. %
We can use the inverse CDF to determine the number of rejections at the 99.9999\% percentile.
For instance, in MATLAB:
\begin{displaymath}
  \sreject = \texttt{icdf(`Negative Binomial', 0.999999, } \szero \texttt{, } \pzero \texttt{)} .
\end{displaymath}
This means that with probability 0.999999,
no more than $\sreject$ nonzeros will be drawn
before $\szero$ zeros are obtained.
So we want to choose the oversampling rate $\rho$ so that
\begin{displaymath}
  \rho \geq (\szero + \sreject) \, \frac{\pzero}{\szero}
  = \frac{\szero + \sreject}{\szero} \, \pzero .
\end{displaymath}

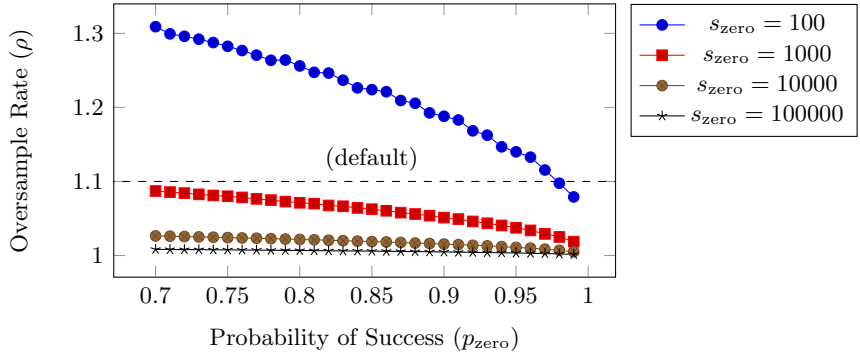
\begin{figure}
  \centering
  \begin{tikzpicture}
    \begin{axis}[
      xlabel=Probability of Success ($\pzero$),
      ylabel=Oversample Rate ($\rho$),
      legend entries={$\szero=100$,$\szero=1000$,$\szero=10000$,$\szero=100000$},
      scale = 0.6,
      label style={font=\footnotesize},
      tick label style={font=\footnotesize},
      legend style={font=\footnotesize},
      width=5 in,
      height=3 in,
      legend pos=outer north east,
      ]
      \addplot table[x=pzero,y=nsamples100] {oversample.dat};
      \addplot table[x=pzero,y=nsamples1000] {oversample.dat};
      \addplot table[x=pzero,y=nsamples10000] {oversample.dat};
      \addplot table[x=pzero,y=nsamples100000] {oversample.dat};
      \draw[dashed] (axis cs:0.6,1.1) -- (axis cs:1.1,1.1)  node[midway,above] {\footnotesize(default)};
    \end{axis}
  \end{tikzpicture}
  \caption{Required oversampling rate for 99.9999\% probability of generating at least $s$ (non-rejected) samples. }
  \label{fig:oversample}
\end{figure}

In \cref{fig:oversample}, we plot the oversample rate that would be needed in different scenarios.
The x-axis is the proportion of nonzeros.
The y-axis is the oversample rate.
We plot four lines corresponding to different values for $\szero$.
We observe two trends:
\begin{enumerate}
\item For fixed $\pzero$, $\rho$ decreases as $\szero$ increases, and
\item For fixed $\szero$, $\rho$ decreases as $\pzero$ increases.
\end{enumerate}
For most real-world examples of sampling zeros from a sparse tensor,
$\pzero \geq 0.99$ because the tensors are extremely sparse.
Additionally, we usually use a sample size of at least $\szero=1000$.
Thus, oversampling by $\rho = 1.1$ should be adequate for most scenarios we expect to encounter.

We could also determine the oversampling rate for any individual problem using this procedure, but  the inverse CDF calculation can be expensive.

\section{Creating Gamma-Distributed Test Problems}\label{sec:creat-gamma-test}%
To create the Gamma-distributed test problem used in \cref{sec:sample-size}, we generate factor matrices whose entries are drawn from the uniform distribution on $(0,1)$:
\begin{displaymath}
  \Ak(i_k,j) \sim \text{uniform}(0,1) \qtext{for all} i_k=1,\dots,n_k, \, j=1,\dots,r, \text{ and } k=1,\dots,d.
\end{displaymath}
Using these factor matrices, we create $\Mtrue$.
The data tensor $\X$ is generated as
\begin{displaymath}
  \xe \sim \text{gamma}(k,\me) \qtext{with} k=1.
\end{displaymath}

\section{Details of Creating Binary Test Problems}%
\label{sec:creating-binary-test}%
We assume an odds link with the data, so the factor matrices must be nonnegative.
The probability of a one is given by $m/(1+m)$ where $m$ corresponds to the odds.
For simplicity in the model and in generating the data tensor,
we assume that factors 1 through $(r-1)$ are relatively sparse (i.e., sparsity specified by $\delta \in (0,1/2)$ and factor $r$ is dense.
The idea here is the last dense component corresponds to noise in the model,
i.e., random but infrequent observations of ones.
Otherwise, the ones have a pattern as dictated by the sparse components.

We specify a density of factor matrix nonzeros and a probability of a one for nonzero values in the resulting model, denoted $\phigh$.
To obtain that probability, the nonzero factor matrix entries should be $\sqrt[d]{\phigh/(1-\phigh)}$.
We modify that slightly by setting the nonzero factor matrix entries to be drawn from a normal distribution with the mean as the target value and a standard deviation of 0.5.
Since only a few entries are nonzero, we can identify all the possible nonzeros corresponding to the first four factors and then compute the exact probability computed by the model and then generate an observation.

For the final dense component, we want the probability of a one, denoted $\plow$, to be relatively low. This means that approximately $\plow$ of the data tensor entries will correspond to this last ``noise'' column. The entries of the factor matrix are set to  $\sqrt[d]{\plow/(1-\plow)}$. We use this value exactly so that we can generate nonzero ``noise'' observations in bulk.

For the test problems in \cref{sec:sample-size}, we use $\delta=0.15$, $\phigh = 0.9$, and $\plow = 0.0025$.
For the test problems in \cref{sec:comp-sampl-strat}, we use $\delta=0.15$, $\phigh = 0.9$, and $\plow = 0.002$.

\section{Cosine Similarity Score}%
\label{sec:meas-recov-true}%
If the true factor matrices are known, we can compute a cosine similarity score between the true and recovered solutions. If the true solution is denoted by $\Ak$ and the estimated solution is $\EstAk$, then the cosine similarity score is
\begin{displaymath}
  \frac{1}{r} \sum_{j=1}^r \prod_{k=1}^d \cos(\Akj,\EstAkj)
\end{displaymath}
where $\pi$ is a permutation that should yield the highest possible similarity. Recall that the cosine of two vectors $\V{a}$ and $\V{b}$ is $\V{a}^{\intercal} \V{b} / (\|\V{a}\|_2 \|\V{b}\|_2)$.
We say that the true solution is recovered if the similarity score is at least 0.9.
Assume that $\M$ holds $\EstAk$ and $\Mtrue$ holds $\Ak$,
the cosine similarity is computed using the Tensor Toolbox for MATLAB \cite{TensorToolbox} via the following command:
\smallskip
\begin{center}\small
\begin{boxedverbatim}
score(M,Mtrue,'lambda_penalty',false)
\end{boxedverbatim}
\end{center}
\smallskip

\section{Individual components of Chicago crime tensor factorization}%
\label{sec:indiv-comp-chic}%
In this appendix, we show the remaining 10 components for the factorization of the Chicago crime tensor discussed in \cref{sec:application-crime} in \cref{fig:comp-1,fig:comp-2,fig:comp-3,fig:comp-4,fig:comp-5,fig:comp-6,fig:comp-7,fig:comp-8,fig:comp-9,fig:comp-10}.
Here we have scaled the date to show the overall weight of the component, and
the other components are normalized.

\morecomppic{1}
\morecomppic{2}
\morecomppic{3}
\morecomppic{4}
\morecomppic{5}
\morecomppic{6}
\morecomppic{7}
\morecomppic{8}
\morecomppic{9}
\morecomppic{10}
%

%
%
%

%
%
%
%

\section{Acknowledgments}%
\label{sec:acknowledgments}%
We thank the referees for their helpful comments on an earlier draft of this paper, including bringing additional references to our attention.
This paper was the result of a broader multi-year collaboration involving several institutions.
We gratefully acknowledge Jed Duersch for preliminary investigations on handling large-scale tensors with GCP and many insightful discussions about the work presented in this paper.
Thanks also to Cliff Anderson-Bergman for stimulating discussions about this project, including providing the information for \cref{sec:determ-overs-rate}.
We thank our colleague Eric Phipps for proposing an idea that ultimately led to semi-stratified sampling.

\bibliographystyle{siamplainmod}


\end{document}